\documentclass[12pt]{article}
\usepackage{epic,amssymb}

\advance \textheight by -1 \topmargin
\advance \textheight by -1 \headheight
\advance \textheight by -1 \headsep
\advance \textheight by -2 \footskip
\vsize = \textheight
\hsize = \textwidth
\newcommand{\Aut}{\mbox{\rm Aut}}
\newcommand{\id}{\mbox{\rm id}}
\newtheorem{theorem}{Theorem}[section]
\newtheorem{proposition}[theorem]{Proposition}
\newtheorem{lemma}[theorem]{Lemma}
\newtheorem{corollary}[theorem]{Corollary}
\newtheorem{definition}[theorem]{Definition}

\newcommand{\bew}{\noindent\underline{Proof.}\ }
\newtheorem{remark}[theorem]{Remark}
\renewcommand{\setminus}{-}
\newcommand{\disj}{\stackrel{.}{\cup}}
\newcommand{\subdir}{\star}
\newcommand{\Z}{{\mathbb{Z}}}
\newcommand{\Q}{{\mathbb{Q}}}
\newcommand{\F}{{\mathbb{F}}}
\newcommand{\N}{{\mathbb{N}}}
\newcommand{\R}{{\mathbb{R}}}

\newcommand{\eb}{\phantom{zzz}\hfill{$\square $}\smallskip}
\renewcommand{\em}{\sf}

\newcommand{\knubbel}{\begin{picture}(0,0)(5,-3)  \put(1,1){\circle*{5}} \end{picture}} 

\title{
 Low dimensional strongly perfect lattices. \\
I: The 12-dimensional case.}

\author{
Gabriele Nebe 
 and Boris Venkov 
\footnote{
The work of Boris Venkov was partially supported by the Swiss National Science Foundation. The work was finished during a visit to the Graduiertenkolleg at
the RWTH Aachen in january 2005. The authors thank the RWTH Aachen for the
invitation.}
} 
 
\begin{document}

 \maketitle

\begin{abstract}
It is shown that the Coxeter-Todd lattice is the unique strongly
perfect lattice in dimension 12.
\end{abstract}

\section{Introduction.}

The notion of perfect lattices came up about 100 years ago in papers 
by Korkine and Zolotarev and  especially Voronoi \cite{Vor} during the
study of dense lattice sphere packings.
If the centers of the spheres in a packing form a lattice $\Lambda $
in Euclidean space $(\R ^n,(,))$
then the  density of the sphere packing is
proportional to the Hermite function of the lattice
$$\gamma (\Lambda ):= \frac{\min ( \Lambda  )}{\det (\Lambda )^{1/n}} $$
where $\min (\Lambda ) :=
\min \{ (\lambda,\lambda ) \mid 0\neq \lambda \in \Lambda  \} $
denotes the square of the minimal distance between distinct lattice points
and $\det (\Lambda )$ is the square of the covolume of $\Lambda $ in 
$\R ^n$.
In any dimension $n$ the Hermite function $\gamma $ 
has a global maximum on the set of $n$-dimensional lattices 
$$\gamma _n := \max \{ \gamma (\Lambda ) \mid \Lambda \subset \R ^n \mbox{ lattice } \} $$
the so called Hermite constant $\gamma _n$.
Exact values for $\gamma _n$ are only known for dimensions 
$n\leq 8$ and, due to a recent work by Cohn and Kumar \cite{Kumar}, 
in dimension $24$.
New upper bounds for $\gamma _n $ are given in the article \cite{Elkies}
by Cohn and Elkies.
The strategy due to Voronoi to find $\gamma _n$ and the densest lattices in 
$\R ^n$ 
is to determine all finitely many local maxima of $\gamma $ on the set of 
similarity classes of $n$-dimensional lattices.
To be a local maximum for $\gamma $, the lattice $\Lambda $ has to be
perfect and eutactic, two conditions on the geometry of the
minimal vectors of $\Lambda $.
For more information the reader is referred to Martinet's book \cite{Martinet}.
Voronoi has developed a remarkable algorithm which permits in principle to
enumerate all (finitely many similarity classes of) perfect lattices
in a given dimension.
But since the number of perfect lattices increases quite rapidely with the
dimension, this seems to be unpracticable in dimensions $\geq 8$.

In \cite{Venkov}, the second author introduced the 
notion of strongly perfect lattices, which are those lattices,
for which the minimal vectors form a 4-design (see Definition \ref{stperf}).
They are perfect and eutactic and hence
local maxima of the Hermite function. Up to rescaling they 
are rational, 
that is the mutual scalar products of lattice vectors lie in $\Q $
(this property is shared by all perfect lattices).
Also most of the famous lattices such as the $E_8$-lattice, the Leech lattice
and the Barnes-Wall lattices are strongly perfect.

There are two general approaches to study and construct strongly perfect
lattices: by modular forms and by invariant theory of finite groups.
The relation with modular forms arises, because the condition that the minimal
vectors of $\Lambda $ form a 4-design means the annulation of certain coefficients in
a theta-series of $\Lambda $ with harmonic coefficients.
In this way one can prove the strong perfectness of many extremal lattices
of small level (see \cite{V} for the even unimodular and \cite{BV} for the
modular case).
By the way this is the only known way to prove the strong perfectness of an 
even unimodular 32-dimensional lattice without roots.
By \cite{King}, there are more than $10^6 $ ot them, an explicit 
classification is not known.

If a rational lattice $\Lambda $ has a big automorphism group 
$G:=\Aut (\Lambda )$ which has no invariant harmonic polynomials of
degree 2 and 4, a condition easily expressed in terms of the character of
$G \leq O(n) $, then $\Lambda $ is strongly perfect.
There are many interesting lattices such as the
Barnes-Wall lattices, the 248-dimensional Thompson-Smith lattice and
others which are strongly perfect by this reason
(see for example \cite{LST}). 

The point of view of lattices is also useful for the study of 
general 4-designs. For example \cite{BMV} treats an infinite 
number of new cases for the classification of tight spherical designs 
by considering the lattice generated by a given design.

The aim of this project, which is a continuation of 
\cite{Venkov}, \cite{dim10}
is to classify all strongly perfect lattices in a given small dimension.
For $n\leq 24$, one only expects few $n$-dimensional strongly perfect lattices
(see \cite[Table 19.1, 19.2, pp. 82,83]{Venkov}) which possibly allow a classification.
For $n\leq 11$ this was done in 
\cite{Venkov}, \cite{dim10}, here we deal with dimension 12. 
Our main theorem is

\vspace{.5cm}

\noindent
{\bf Theorem.} (see Theorem \ref{main12})
{\it 
The strongly perfect lattices in dimension $12$ are similar to the
Coxeter-Todd lattice $CT$.}

\vspace{.5cm}

We prove this theorem by eliminating rather easily all possibilities for
the Kissing number $2s = |\Lambda _{\min } |$ of a strongly perfect lattice
in dimension 12, except for $s=252$ and $s=378$.
To eliminate $s=252$ is more difficult and involves the construction
of higher-dimensional lattices up to dimension 18.
In the last case, $s=378$, the lattice $\Lambda $ has the same 
Kissing number as the Coxeter-Todd lattice.
One special property of the Coxeter-Todd lattice, is that it is 
$3$-modular and has a natural complex structure as a unimodular
hermitian lattice over $\Z [\frac{1+\sqrt{-3}}{2}]$.
Its shortest vectors form the root system of the complex reflection
group number 34 $(\cong 6.U_4(3).2)$ in the Shephard-Todd classification
\cite{ShephardTodd}.
During the identification of the putative strongly perfect lattice $\Lambda $
with $s=378$ we construct step by step parts of this rich structure
which finally allows us to show that $\Lambda $ is 
similar to the Coxeter-Todd lattice.

\section{Some general equations.}{\label{allg}}

\subsection{General notation.}{\label{not}}

For a lattice $\Lambda $ in $n$-dimensional Euclidean space we denote by
$\Lambda ^*$ its dual lattice and 
by
$$\Lambda _a :=\{ \lambda \in \Lambda \mid (\lambda , \lambda ) = a \} $$
the vectors of square length $a$. 
In particular, 
$\Lambda _{\min }$ is the set of minimal vectors in $\Lambda $, its
cardinality is known as the {\em Kissing number} of the lattice $\Lambda$.
There are general bounds on the cardinality of an 
antipodal spherical code and hence on the 
 Kissing number of an $n$-dimensional lattice. 
For $n=12$ the bound is
 $ |\Lambda _{\min } |\leq 2\cdot 614 $
 (see \cite[Table 1]{Anstreicher}).

\subsection{Designs and strongly perfect lattices}

Let $(\R ^n, (,))$ be the Euclidean space of dimension $n$.
For  $m\in \R $, $m>0$ denote by 
$$ S^{n-1}(m) := \{
y\in \R^n \mid (y,y)=m \}$$ 
the $(n-1)$-dimensional sphere of radius $\sqrt{m}$.

\begin{definition}
A finite nonempty set 
$ X \subset S^{n-1}(m)$ is called a {\em spherical $t$-design},
if 
$$\frac{1}{|X|} \sum _{x \in X} f(x) = \int _{S^{n-1}(m)} f(x) d\mu (x) $$
for all polynomials of $f\in \R [x_1,\ldots , x_n ]$ of degree $\leq t$,
where $\mu $ is the  $O(n)$-invariant measure on the sphere, normalised 
such that 
$\int _{S^{n-1}(m)} 1 d\mu (x) = 1 $.
\end{definition}

Since the condition is trivially satisfied for constant polynomials $f$,
and the harmonic polynomials generate the orthogonal complement 
$\langle 1 \rangle ^{\perp } $ with respect to the
$O(n)$-invariant scalar product 
$<f,g>:= \int _{S^{n-1}(m)} f(x) g(x) d \mu (x)  $ on $\R [x_1,\ldots , x_n]$,
it is equivalent to ask that 
$$\sum _{x\in X } f(x) = 0 $$ 
for all harmonic polynomials $f $ of degree $\leq t$.

\begin{definition}{\label{stperf}}
A lattice $\Lambda  \subset \R^n$ is called 
{\em strongly perfect}, if its minimal vectors 
$\Lambda _{\min }  $ form a spherical $4$-design.
\end{definition}

The most important property of strongly perfect lattices is that
they provide interesting examples for local maxima of the Hermite 
function (see \cite{Venkov}). 

Let $\Lambda $ be a strongly perfect lattice of dimension $n$,
$m:= \min (\Lambda )$ and choose
$X\subset \Lambda _{m}$ such that $X\cup -X = \Lambda _{m}$
and $X\cap -X = \emptyset$.
Put $s:= |X| $.

By \cite{Venkov} the condition that $\sum _{x\in X} f(x) = 0$ for all
harmonic polynomial of degree 2 and 4 may be reformulated 
to the condition that for all $\alpha \in \R^n$
$$ (D4) (\alpha ): \ \ \sum _{x\in X} (x,\alpha )^4 = \frac{3sm^2}{n(n+2)} (\alpha , \alpha )^2.$$
Applying the Laplace operator $\sum _{i=1}^n \frac{\partial ^2}{\partial \alpha _i^2 }$ to $(D4)(\alpha )$  one obtains
$$ (D2) (\alpha ): \ \ \sum _{x\in X} (x,\alpha )^2 = \frac{sm}{n} (\alpha , \alpha ) .$$

Note that it is enough to assume that  there are constants $c_t$ ($t=1,2$) such
that 
$$\sum _{x\in X} (x,\alpha )^{2t} = c_t (\alpha , \alpha ) ^t \mbox{ for all }
\alpha \in \R^n .$$ 
These constants $c_t$ are then uniquely determined as one sees 
applying $t$ times the Laplace operator with respect to $\alpha $.
For later use, we also remark that the condition $(D2)$ is equivalent to
the $2$-design property of $X\cup -X$.

Substituting $\alpha := \xi _1 \alpha _1 + \xi _2 \alpha _2$ 
in  $(D2)$  and comparing coefficients, one finds
$$(D11)(\alpha _1,\alpha _2) \ \ \sum _{x\in X} (x,\alpha _1 ) (x,\alpha _2)= \frac{sm}{n} (\alpha _1, \alpha _2 ) \mbox{ for all } \alpha _1 ,\alpha _2 \in \R^n.$$
Writing $\alpha $ as a linear combination of 4 vectors, $(D4)$ implies that
for  all $\alpha _1,\ldots , \alpha _4 \in \R^n$
$$(D1111)
\ \sum _{x\in X} (x,\alpha _1 ) (x,\alpha _2) (x,\alpha _3) (x,\alpha _4) 
= $$
$$ = \frac{sm^2}{n(n+2)} (
(\alpha _1, \alpha _2 )(\alpha _3 , \alpha _4) +
(\alpha _1, \alpha _3 )(\alpha _2 , \alpha _4) +
(\alpha _1, \alpha _4 )(\alpha _2 , \alpha _3))
$$
In particular
$$(D13)(\alpha _1,\alpha _2): \ \ \sum _{x\in X} (x,\alpha _1 ) (x,\alpha _2)^3= \frac{3sm^2}{n(n+2)} (\alpha _1, \alpha _2 )(\alpha _2 , \alpha _2)$$
$$(D22)(\alpha _1,\alpha _2): \ \ \sum _{x\in X} (x,\alpha _1 )^2 (x,\alpha _2)^2= \frac{sm^2}{n(n+2)} 
(2(\alpha _1, \alpha _2 )^2+(\alpha _1 , \alpha _1)(\alpha _2 , \alpha _2))$$
$(D13)-(D11)(\gamma ,\alpha )$ reads as 
$$ \sum _{x\in X} (x,\gamma ) (x,\alpha )
((x,\alpha )^2-1)
= (\alpha, \gamma  )\frac{sm}{n}
( \frac{3m}{n+2} (\alpha  , \alpha ) - 1 ) \mbox{ for all } \alpha , \gamma \in \R^n ~.$$
Note that for $\alpha \in \Lambda ^*$ and $x\in \Lambda $, the product 
$(x , \alpha ) ((x,\alpha )^2-1) $ is divisible by $6$ and hence 
$\frac{1}{6} ((D13)-(D11))(\gamma ,\alpha ) \in \Z $  yields that $$ 
\frac{1}{6} \sum _{x\in X} (x,\gamma ) (x,\alpha )
((x,\alpha )^2-1)
= (\alpha, \gamma  )\frac{sm}{6n}
( \frac{3m}{n+2} (\alpha  , \alpha ) - 1 )  \in \Z \ \mbox{ for all } \alpha , \gamma \in \Lambda ^* .$$
Putting $\alpha = \gamma $ in $(D13) - (D11)$ we obtain 
$\frac{1}{12} ((D4)-(D2))(\alpha )  \in \Z $ hence $$
\frac{1}{12} \sum _{x\in X} (x,\alpha )^2
((x,\alpha )^2-1) = \frac{sm}{12n} (\alpha, \alpha  )
( \frac{3m}{n+2} (\alpha  , \alpha ) - 1 )  \in \Z \ \mbox{ for all } \alpha  \in \Lambda ^* $$
since 
$(x , \alpha )^2 ((x,\alpha )^2-1) $ is divisible by $12$ if $(x,\alpha )\in \Z$.

\begin{lemma}{\label{linkomb}}(see \cite[Lemma 2.1]{dim10})
Let $\alpha \in \R ^n$ be 
such that $(x , \alpha ) \in \{ 0,\pm 1 ,\pm 2 \} $ for all $x\in X$.
Let $N_2(\alpha ) := \{ x\in X\cup -X  \mid (x,\alpha ) = 2 \}$ and
put $$c:=\frac{sm}{6n}(\frac{3m}{n+2}(\alpha ,\alpha )  -1 ) .$$
Then 
$|N_2(\alpha ) | = c (\alpha , \alpha ) /2 $ and
$$\sum _{x\in N_2(\alpha )} x = c \alpha . $$
\end{lemma}

Lemma \ref{linkomb} will be often applied to $\alpha \in \Lambda ^*$.
Rescale $\Lambda $ such that $\min (\Lambda ) = m = 1$
and let $r:=\min (\Lambda ^*)$.
Since $\gamma(\Lambda ) \gamma (\Lambda ^*) = \min(\Lambda )
\min (\Lambda ^*) \leq \gamma _n ^2$, we 
get $r\leq \gamma _n^2$ and for $\alpha \in \Lambda ^*_r$ we
have $(\alpha , x)^2 \leq r$ for all $x\in \Lambda _1$.
Hence if $r<9$ then $(\alpha , x ) \in \{ 0,\pm 1 ,\pm 2\}$ for all
$x\in X$ and Lemma \ref{linkomb} may be applied.

The next lemma yields  good bounds on $n_2(\alpha )$.

\begin{lemma}{\label{boundn2}}
Let $m:= \min (\Lambda ) $ and choose $ \alpha \in \Lambda ^*_r$.
If $r\cdot m < 8$, then
 $$|N_2(\alpha ) |\leq \frac{r m}{8-r m} .$$
\end{lemma}

\bew
Since $(x,x)(\alpha , \alpha ) < 8$ for all $x\in \Lambda _m$, the
scalar product $|(x,\alpha )| < 3$. Hence $\alpha $ satisfies the
conditions of Lemma \ref{linkomb}.
Let $N_2(\alpha ) = \{ x_1,\ldots , x_k \}$ and $c =\frac{2k}{r}$
be the constant from Lemma \ref{linkomb}.
Then $(x_i,x_i ) = m$ and
$(x_i,x_j) \leq \frac{m}{2}$ because the $x_i$ are minimal vectors in $\Lambda $.
Hence 
$$ \frac{4k}{r}  = (x_1,c\alpha ) = (x_1,x_1) + \sum _{i=2}^k (x_1,x_i ) 
\leq m + \frac{m(k-1)}{2} = \frac{m(k+1)}{2} $$
which yields that 
$k=|N_2(\alpha )| \leq \frac{rm}{8-rm} .$
\eb

\begin{lemma}{\label{min}}(\cite[Th\'eor\`eme 10.4]{Venkov})
Let $L$  be a strongly perfect lattice of dimension $n$.
Then 
$$\gamma(L) \gamma (L^*) =  \min(L) \min(L^*) \geq \frac{n+2}{3} .$$
A strongly perfect lattice $L$ where equality holds is called 
{\em of minimal type}.
\end{lemma}

As an application this lemma allows to show that $|N_2(\alpha ) | \neq 1$.

\begin{lemma}
Choose $ \alpha \in \Lambda ^*_r$ that
satisfies the conditions of Lemma \ref{linkomb}.
If $n\geq 11$ then $|N_2(\alpha )| \neq 1 $.
\end{lemma}

\bew
Let $m:=\min(\Lambda )$ and assume that $N_2(\alpha ) = \{ x \}$.
Then $\alpha = c x $ for some constant $c$.
Taking scalar products with $x$ yields $c=\frac{2}{m}$
and hence  $(\alpha , \alpha ) = \frac{4}{m}$.
Therefore 
$$\frac{n+2}{3} \leq \min  (\Lambda ) \min (\Lambda ^*)  \leq 4 $$
which implies that $n\leq 10$.
\eb

{\bf Notation:} 
For a set $X$ as above (usually clear from the context)
and $\alpha \in \R ^n$ we let
$$N_i(\alpha ) = \{ x\in X \mid (\alpha , x) = \pm i \} .$$

The last lemma is quite useful in the investigation of lattices
of minimal type:

\begin{lemma}\label{2design}
Let $X \disj -X \subset S^{n-1}(m) $ be a spherical $4$-design and let $\alpha \in \R^n$ 
be such that $(x,\alpha ) \in \{ 0,\pm 1 \}$ for all $x\in X$.
Let $M:=N_1 (\alpha ) $ and
let $\pi : \R ^n \to \langle \alpha \rangle ^{\perp }$ be the orthogonal projection
onto the orthogonal complement of $\alpha $.
Then $\pi (M) \disj -\pi (M)  \subset S^{n-2}(m' )$ is a spherical $2$-design.
\end{lemma}

\bew
Let $v\in \langle \alpha \rangle ^{\perp }$, i.e. $(v,\alpha ) = 0$.
Then by $(D22) (\alpha , v) $ we find that 
$$\sum _{x\in M } (x,v)^2 = \sum _{x\in \pi (M) } (x,v)^2 = c (v,v) $$
for some constant $c$ not depending on $v$.
Therefore $\pi (M) \cup - \pi(M) $ is a $2$-design.
\eb

\subsection{Some general facts on lattices.}

The next two lemmas about indices of sublattices are used quite often
in the argumentation below.
Since we are dealing with norms modulo some prime number $p$, we may
pass to the localization
$\Z _p := ( \Z \setminus p\Z )^{-1} \Z \subset \Q $ of
$\Z $ at $p$.

\begin{lemma}\label{index4}
Let $\Gamma $ be a  $\Z _2$-lattice such that $(\gamma ,\gamma ) \in \Z _2 $ for 
all $\gamma \in \Gamma $.
Let $\Gamma ^{(e)} := \{ \alpha \in \Gamma \mid (\alpha,\alpha ) \in 2 \Z _2 \}$.
If $\Gamma ^{(e)} $ is a sublattice of $\Gamma $, then 
$[\Gamma : \Gamma ^{(e)}] \in \{ 1,2,4 \}$.
\end{lemma}

\bew
Clearly $2\Gamma \subset \Gamma ^{(e)}$ hence $\Gamma / \Gamma ^{(e)} $ is a
vector space over $\F _2 $.
Moreover $(\alpha  ,\beta ) \in \frac{1}{2} \Z _2 $ for all $\alpha,\beta \in \Gamma $ since the norms in $\Gamma $ are integers.
Let $\alpha , \beta, \gamma \in \Gamma \setminus \Gamma ^{(e)} $.
Then $(\alpha + \beta , \alpha + \beta ) =
(\alpha ,\alpha ) + (\beta , \beta ) + 2 (\alpha ,\beta ) \in
2(\alpha , \beta ) + 2 \Z _2$ since both
$(\alpha,\alpha ) $ and $(\beta ,\beta )$ are odd.
In particular, if $(\alpha , \beta ) \in \Z _2$, then $\alpha + \beta \in \Gamma ^{(e)}$.
Since one of $(\alpha ,\gamma) $, $(\beta ,\gamma )$ or 
$(\alpha + \beta , \gamma ) $ is integral, the three classes
$\alpha +\Gamma ^{(e)}, \beta + \Gamma ^{(e)}, \gamma + \Gamma ^{(e)} $ are
$\F _2$-linearly dependent and $\dim _{\F _2}(\Gamma /\Gamma ^{(e)} )\leq 2$.
\eb

\begin{lemma}\label{index3}
Let $\Gamma $ be a $\Z _3$-lattice such that $(\gamma ,\gamma ) \in \Z _3 $ for 
all $\gamma \in \Gamma $.
Let $\Gamma ^{(t)} := \{ \alpha \in \Gamma \mid (\alpha,\alpha ) \in 3 \Z _3 \}$.
Assume that 
$$ (\alpha ,\beta )^2 - (\alpha , \alpha ) (\beta , \beta)  \in 3\Z_3 \mbox{ for
all } \alpha , \beta \in \Gamma .$$
Then 
$\Gamma ^{(t)} $ is a sublattice of $\Gamma $ and
$[\Gamma : \Gamma ^{(t)}] \in \{ 1,3 \}$.
\end{lemma}

\bew
 By our assumption $(\alpha  ,\beta ) \in 3 \Z $ for all $\alpha,\beta \in \Gamma $  with
$\alpha \in \Gamma ^{(t)}$.
Therefore $\Gamma ^{(t)}$ is a sublattice of $\Gamma $.
Clearly $3\Gamma \subset \Gamma ^{(t)}$ hence $\Gamma / \Gamma ^{(t)} $ is a
vector space over $\F _3 $.
Let $\alpha , \beta\in \Gamma \setminus \Gamma ^{(t)} $.
Then $(\alpha , \alpha) (\beta ,\beta ) \equiv (\alpha , \beta )^2 \pmod{3}$
implies that $(\alpha ,\alpha ) \equiv (\beta , \beta )\pmod{3} $ and
one of $\alpha \pm \beta \in \Gamma ^{(t)}$.
Hence $[\Gamma  : \Gamma ^{(t)} ] \leq 3$.
\eb

We will also meet 
families of vectors  
$E:=\{ v_1,\ldots , v_k \}$ of equal norm, say $(v_i,v_i) = 1$ and 
non positive scalar products, i.e. 
$(v_i,v_j) \leq 0$ for all $i\neq j$
with $\sum _{i=1}^k v_i = 0$.
Such a system is called {\em decomposable}, if $E= E' \disj E''$ such that
$$\sum _{v\in E'} v = \sum _{v\in E''} v = 0 .$$
It is classical (and easy) that for any indecomposable system $E$, the
relation $\sum _{ v\in E} v = 0$ is the only relation between the 
vectors in $E$.
In particular 
$$k = |E| = \dim \langle E \rangle + 1 \mbox{ for any indecomposable } E. $$
An arbitrary system $E$ can be written as an orthogonal union 
$E= \disj _{i=1}^t E_i $ with $E_i \perp E_j$, $E_i$ indecomposable.

\begin{lemma}{\label{neg}}
If $t$ is the number of indecomposable components of $E$, then
$|E| = \dim \langle E \rangle +  t $
\end{lemma}

\bew
For $t=1$ this is clear by the argumentation above.
For general $E$, all three functions 
$|E|$, $\dim \langle E \rangle $, $t(E)$ are additive for 
orthogonal union.
\eb

Components of dimension 1 are just pairs $\{ x, -x \}$.
Components of dimension 2 are of the form $\{ x_1,x_2,x_3 \}$ with 
$x_1+x_2+x_3 = 0$ and $(x_i,x_i) = 1$, $(x_i,x_j) = -\frac{1}{2} $ for
all $i\neq j$.
They generate a root system $A_2$.
For more information see \cite{CoSVI}.

\section{Dimension 12}

In this section we prove our main theorem:

\begin{theorem}\label{main12}
Let $\Lambda $ be a strongly perfect lattice of dimension $12$.
Then $\Lambda \cong CT$ is similar to the Coxeter-Todd lattice.
\end{theorem}

The lattice $CT$, usually denoted by $K_{12}$,
is an extremal 3-modular lattice in the
sense of \cite{Queb}, which means that 
$CT$ is similar to its dual lattice.
Rescaled to minimum $4$, the determinant of $CT$ is $3^6$.
The Kissing number is $|CT_{4} | = 2 \cdot 378 =756$.
The Coxeter-Todd lattice is the densest known lattice in dimension 12.

\begin{remark}\label{CTmin}
No proper overlattice of $CT$ has minimum $4$.
\end{remark}

\bew
Let $\Gamma > CT  $ be a proper overlattice of $CT$ with 
$\min (\Gamma ) = 4$. Then $\det (\Gamma ) = \det (CT) \cdot [\Gamma : CT]^{-2} 
\leq 3^6 \cdot 2^{-2} $.
Hence the Hermite function 
$$\gamma (\Gamma ) = \frac{4}{\sqrt[12]{\det(\Gamma )}} \geq 
2.59 > \gamma _{12} $$
which is a contradiction.
\eb

Therefore it is enough to show that any strongly perfect lattice
of dimension 12, that is generated by its minimal vectors, is similar
to $CT$.

\subsection{Kissing numbers.}

Let $\Lambda $ be a strongly perfect lattice in dimension 12, rescaled
such that $\min (\Lambda ) = 1$. 
Then by Lemma \ref{min} and the Cohn-Elkies bound that $\gamma _{12} \leq 2.522 $
we find 
$$ \frac{14}{3} \leq r:=\min(\Lambda ^* ) \leq \gamma _{12}^2 \leq 6.37 .$$
Hence $\alpha \in \Lambda ^* _r$ satisfies the hypothesis of Lemma \ref{linkomb}
and Lemma \ref{boundn2} yields 
$$n_2(\alpha ) \leq \frac{r}{8-r} \leq \frac{\gamma _{12}^2}{8-\gamma _{12}^2}
\leq 3.88 < 4 .$$

If $d:=\det (\Lambda ^* ) = \frac{1}{\det (\Lambda )}$ denotes the
determinant of $\Lambda ^* $, then 
$d\leq \gamma _{12}^{12} \leq 66212.7 $.

Let $s:= |X|$ where $X\disj -X = \Lambda _1$ be half the 
Kissing number of $\Lambda $.
Then $\frac{12\cdot 13}{2} = 78 \leq s \leq 614$ and 
$r$ is a rational solution of 
$$ n_2(\alpha ) = \frac{sr}{12\cdot 12\cdot 14} (3r-14) .$$
Going through all possibilities by a computer we find:

\begin{proposition}{\label{poss}}
With the notation above, one of the following holds:
\begin{itemize}
\item[(a)] $s=168= 2^3\cdot 3\cdot 7 $ and $r=6$.
\item[(b)] $s=252 = 2^2\cdot 3^2 \cdot 7$ and $ r=6 $.
\item[(c)] $s=378 = 2 \cdot 3^3 \cdot 7$ and $r=6$.
\item[(d)] $r=\frac{14}{3}$ and $\Lambda $ is of minimal type.
\end{itemize}
\end{proposition}

\subsection{The case $s=168$, $r=6$.}

\begin{proposition}\label{b}
There is no strongly perfect lattice satisfying Proposition \ref{poss} (a).
\end{proposition}

\bew
Let $\Lambda $ be such a strongly perfect lattice with
$\min (\Lambda ^* ) = r = 6$ and $s(\Lambda ) = s = 168 $.
Then for all $\alpha \in \Lambda ^* $
$$
\begin{array}{lll} \sum _{x\in X }(x,\alpha )^2  & = 14 (\alpha , \alpha ) & \in \Z \\
\sum _{x\in X }(x,\alpha )^4 & = 3 (\alpha , \alpha )^2& \in \Z  
\end{array}
$$
hence all norms in $\Lambda ^* $ are integers.
The equalities 
$\frac{1}{6} ((D13) - (D11))$ and $\frac{1}{12} ((D2)-(D4)) $
yield that 
$$
\begin{array}{ll}
\frac{1}{6} (\alpha _1,\alpha_2) (3 (\alpha _2,\alpha _2) -14) & \in \Z \\
\frac{1}{12} (\alpha _1,\alpha_1) (3 (\alpha _1,\alpha _1) -14) & \in \Z 
\end{array}
$$
for all $\alpha _1,\alpha _2\in \Lambda ^* $.
Therefore $\Lambda ^* $ is an even lattice with $3\mid (\alpha _1 ,\alpha _2)$
for all $\alpha _1,\alpha _2 \in \Lambda ^* $ and hence 
$\Gamma := \frac{1}{\sqrt{3}} \Lambda ^* $ is an even lattice with
$\min (\Gamma ) = 2 $ and $\min (\Gamma ^* ) = 3$ which is a contradiction.
\eb

\subsection{Lattices of minimal type.}

In this section we prove the following 

\begin{theorem}
Let $\Lambda $ be a strongly perfect lattice of dimension $12$ and of minimal 
type, i.e. 
$\min (\Lambda ) \min (\Lambda ^*) = \frac{n+2}{3} = \frac{14}{3} $
and let $s:= \frac{1}{2} | \Lambda _{\min } |$.
Then $s = 378 $
or $s= 252$.
\end{theorem}

Rescale $\Lambda $, such that $\min (\Lambda) = 1$.
Then $\min (\Lambda ^*) = \frac{14}{3} $.
Applying $(D2)$ to $\alpha \in (\Lambda ^*)_{14/3} $ yields that 
$$\frac{14 s}{12\cdot 3 } = \frac{7s}{18} \in \Z .$$
Hence $$s=18 s_1 \leq 614 \mbox{ for some  } s_1 \in \N, s_1 \leq 34 .$$

\begin{lemma}
$7$ divides $s_1$.
\end{lemma}

\bew
Assume that $7$ does not divide $s_1$ and choose $\alpha \in \Lambda ^*$.
The $(\alpha , \alpha ) = \frac{p}{q}$ for some coprime integers $p,q$.
Since  $(D4)(\alpha )  = 3^2\cdot 2^{-2} \cdot 7^{-1} s_1 \frac{p^2}{q^2}$
is integral, this implies that $p = 7 p_1$ is
divisible by $7$ and $q^2$ divides $3^2 s_1$. 
Moreover $((D4) - (D2))(\alpha )$ is divisible by 12 which yields  that 
$$ (\star ) \ \ \frac{s_1}{2^4 q^2} 7p_1 (3p_1-2q)   \in \Z $$
If $q$ is even then $p_1$ is odd and $2^6 = 64 \mid 2^4 q^2 \mid s_1$ which contradicts
the fact that $s_1 \leq 39$.
Since $q^2$ divides $3^2 s_1$, the only possibilities for $q$ are
$$q = 1 , 3,3^2,3 \cdot 5  .$$

\knubbel
First assume that there is $\alpha \in \Lambda ^*$ with $(\alpha , \alpha ) = \frac{7p_1}{15}$.

Then $s_1 = 5^2$ is odd and hence by $(\star )$ the norms of the elements in 
$\Lambda ^*$ lie in $\frac{14}{15} \Z $.
Let $$\Gamma := \sqrt{\frac{15}{14}} \Lambda ^* .$$
Then all norms in $\Gamma $ are integers.
For the scalar products we apply $(D22)$ to $\alpha , \beta \in \Gamma $ to find 
that 
$$(\star \star) \ \ \frac{7}{3} (2 (\alpha , \beta )^2 + (\alpha , \alpha ) (\beta ,\beta )) \in \Z .$$
Therefore $\Gamma $ is an integral lattice.
Let $$
\begin{array}{rl}
\Gamma ^{(e)} :=  & \{ \gamma  \in \Gamma \mid (\gamma ,\gamma ) \in 2 \Z \} \\
\Gamma ^{(t)} :=  & \{ \gamma  \in \Gamma \mid (\gamma ,\gamma ) \in 3 \Z \} \\
\Gamma ':=\Gamma ^{(e)} \cap \Gamma ^{(t)} :=  & \{ \gamma  \in \Gamma \mid (\gamma ,\gamma ) \in 6 \Z \} 
\end{array}
$$
Then $\Gamma ^{(e)}$ is the even sublattice of $\Gamma $ 
of index $2$.
Moreover 
$((D4)- (D2) ) (\alpha )$ yields that 
$ 
 \frac{7}{12} (\alpha , \alpha ) ( (\alpha ,\alpha ) - 5 ) \in \Z $
 for all $ \alpha \in \Gamma $.
 In particular 
the norms in $\Gamma ^{(e)}$ are divisible by
$4$ and hence $\frac{1}{\sqrt{2}} \Gamma ^{(e)}$ is still even.
By $(\star \star) $ and Lemma \ref{index3} 
$\Gamma ^{(t)}$ is a sublattice of 
index $\leq 3$ in $\Gamma $.
In particular $[\Gamma : \Gamma '] \leq 6$ and $$\Delta := \frac{1}{\sqrt{6}} \Gamma '$$ is even.
Since $\min (\Gamma ^*) \det (\Gamma )^{1/12} \leq \gamma _{12} \leq 2.522$ we get
$$\det (\Gamma ) \leq (2.522 \frac{15}{14} )^{12} < 151530.4 $$
and therefore 
$$\det (\Delta ) \leq \frac{1}{6^{12}} 6^2 \det (\Gamma ) \leq \frac{151530}{6^{10}} \leq 0.026 $$
which is a contradiction since $\Delta $ is integral.

\knubbel
Now assume that there is $\alpha \in \Lambda ^*$ with $(\alpha , \alpha ) = \frac{7p_1}{9}$

Then $s_1 = 3^3$ is odd and hence by $(\star )$ the norms of the elements in
$\Lambda ^*$ lie in $\frac{14}{9} \Z $.
Let $$\Gamma := \sqrt{\frac{9}{14}} \Lambda ^* .$$
Then $\min (\Gamma ) = 3$, $\min (\Gamma ^*) = \frac{14}{9}$ and all norms in $\Gamma $ are integers.
For the scalar products we apply $(D22)$ to $\alpha , \beta \in \Gamma $ to find
that
$$7 (2 (\alpha , \beta )^2 + (\alpha , \alpha ) (\beta ,\beta )) \in \Z .$$
Therefore $\Gamma $ is an integral lattice.

In the new scaling the equation
$((D4) - (D2)) (\alpha )$ yields that
$$ 
 \frac{7}{4}   (\alpha , \alpha ) ( (\alpha ,\alpha ) - 3 ) \in \Z  \mbox{ for all } \alpha \in \Gamma
.
$$

Let $$\Gamma ^{(e)} := \{ \gamma \in  \Gamma \mid (\gamma , \gamma ) \in 2 \Z  \} $$
be the even sublattice of $\Gamma $.
Then the norms in $\Gamma ^{(e)}$ are divisible by $4$  and hence
$\frac{1}{\sqrt{2}} \Gamma ^{(e)} $ is an integral lattice.
Therefore
$$1024 = 2^{10} \leq \frac{\det(\Gamma ^{(e)})}{4} = \det(\Gamma )\leq (2.522 \frac{9}{14} )^{12} \leq 330 $$
which is a contradiction.

\knubbel
The remaining case is that all norms of elements in $\Lambda ^*$ lie in $\frac{7}{3} \Z $.

Let $\Gamma := \sqrt{\frac{3}{7}} \Lambda ^*$.
Then $\min (\Gamma ) = 2$, $\min (\Gamma ^*) = \frac{7}{3} > 2$ and all norms in $\Gamma $ are integers.
But $\Gamma $ is not an integral lattice since $\min (\Gamma ) < \min (\Gamma ^*)$.
In particular, there is $\alpha \in \Gamma $ with $(\alpha ,\alpha ) \in 1 + 2\Z $.

In the new scaling the equation
$\frac{1}{12} ((D4) - (D2)) (\alpha )$ and $\frac{1}{6} ((D13)- (D11))(\alpha , \beta )$ yield that
$$
\begin{array}{l} 
 2^{-4}\cdot 3^{-1} \cdot 7 s_1 (\alpha , \alpha ) ( (\alpha ,\alpha ) - 2 ) \in \Z \mbox{ for all } \alpha \in \Gamma \\
2^{-3}\cdot 3^{-1} \cdot 7 s_1 (\alpha , \beta ) ( (\alpha ,\alpha ) - 2 ) \in \Z  \mbox{ for all } \alpha, \beta \in \Gamma 
.
\end{array}
$$
In particular 
$$2^4 \mid s_1 \in \{ 16, 32 \} \mbox{ and } 
3 \mbox{ divides }  (\alpha ,\beta ) ((\alpha , \alpha ) - 2 ) \mbox{ for all } \alpha ,\beta \in \Gamma .$$
As above we find that 
$$\Gamma ^{(t)} := \{ \gamma \in \Gamma \mid (\gamma , \gamma ) \in 3 \Z \} $$
is a sublattice of $\Gamma $ of index $3$ such that 
$\sqrt{\frac{2}{3}} \Gamma ^{(t)}$ is even.
Therefore
$$14 \leq \frac{3^{10}}{2^{12}} \leq \frac{\det(\Gamma ^{(t)})}{9} = \det(\Gamma )\leq (2.522 \frac{3}{7} )^{12} \leq 2.6  $$
which is a contradiction.
\eb

We therefore have $s = 2\cdot 3^2 \cdot 7 \cdot s_2 $ with 
$s_2 \in \{ 1,2,3,4 \}$.
To obtain the theorem it remains to show that 
$s_2=3$ or $s_2=2$.
We keep the scaling such that $\min (\Lambda ) = 1$.

\begin{lemma}\label{gamma1}
Assume that $s_2 \neq 3$. 
Then 
$$\Gamma ^{(t)} := \{ \gamma \in \Lambda ^* \mid (\gamma , \gamma ) \in \Z \} $$
is an even sublattice of $\Lambda ^*$
of index $[\Lambda ^* : \Gamma ^{(t)}] = 3 $.
\end{lemma}

\bew
For $\alpha \in \Lambda ^*$ we write $(\alpha , \alpha ) = \frac{p}{q}$ with
coprime integers $p,q$. Then $\frac{1}{12} ((D4) - (D2))(\alpha ) \in \Z$ 
reads as
$$
(\star \star ) \ \ \frac{s_2}{2^4q^2} p(3p-14q) \in \Z .$$
Since $2^4$ does not divide $s_2$, we find that $p$ is even and 
$q\in \{ 1,3 \}$.

Moreover 
$\frac{1}{6} ((D13) - (D11))(\alpha , \beta )$ yields that 
$$(\star _2) \ \ \frac{1}{2^3} s_2 (\alpha , \beta ) (3 (\alpha , \alpha ) -14 )  \in \Z 
\mbox{ for all } \alpha , \beta \in \Lambda ^* .$$
Hence if $s_2 \neq 3$, then 
$(\alpha , \beta ) \in \Z $ for all $\alpha \in \Lambda ^*$ for
which $(\alpha,\alpha ) \in \Z$, therefore
$$\Gamma ^{(t)} := \{ \gamma \in \Lambda ^* \mid (\gamma , \gamma ) \in \Z \} $$
is an even sublattice of $\Lambda ^*$.
The index of $\Gamma ^{(t)}$ in $\Lambda ^*$ is $3$, since 
for $\alpha ,\beta \in \Lambda ^* \setminus \Gamma ^{(t)}$, 
$(D22)(\alpha , \beta )$ implies that  
$$\frac{3s_2}{4} (2 (\alpha , \beta )^2 + (\alpha , \alpha ) (\beta , \beta )) \in \Z.$$
Rescaling the bilinear form with 3, this allows to apply Lemma 
\ref{index3} to see that $[\Lambda ^* : \Gamma ^{(t)} ] = 3 $.
\eb

\begin{lemma}{\label{s2=1,5}}
$s_2 \neq 1$.
\end{lemma}

\bew
Assume that $s_2 = 1 $.
The norms of the elements of the even lattice $\Gamma ^{(t)}$ in
Lemma \ref{gamma1} satisfy $(\star \star )$ with $q=1$ and $s_2 = 1$.
Moreover $\min (\Gamma ^{(t)}) \geq \frac{14}{3} $.
Since $p=6$ does not satisfy $(\star \star )$ we find
$\min (\Gamma ^{(t)}) \geq 8$.
Since $d:= \det (\Lambda ^* ) \leq \gamma _{12}^{12} \leq 66212.7 $ we get
$$\gamma (\Gamma ^{(t)}) \geq \frac{8}{ (9 d)^{1/12}} \geq \frac{8}{(9\cdot 66212.7)^{1/12}}\geq 2.64$$
contradicting that $\gamma _{12} \leq 2.522 $.
\eb

\begin{lemma}
$s_2 \neq 4$.
\end{lemma}

\bew
Assume that $s_2 = 4$.
Then $\Gamma ^{(t)} $ is an even lattice of minimum $\geq \frac{14}{3} $.
Hence $\min (\Gamma ^{(t)} ) \geq 6$.
If we show that $\Gamma^{(t)} $ does not contain vectors of square length 6, then
we obtain a contradiction as in Lemma \ref{s2=1,5}.
So let $\gamma \in \Gamma ^{(t)}$ be an element of norm $(\gamma ,\gamma ) = 6$.
Then for all $x\in X$, $(x,\gamma ) \in \{ 0,\pm 1,\pm 2\}$ and hence by
Lemma \ref{linkomb} $n_2(\gamma ) = 6 > \frac{6}{8-6} =3$
contradicting the bound of Lemma \ref{boundn2}.
\eb

\subsection{R\'esum\'e}

Let us summarise what we have shown:

\begin{theorem}
Let $\Lambda $ be a strongly perfect lattice of dimension $12$ and 
$s:= \frac{1}{2} | \Lambda _{\min } |$ be half the kissing number of $\Lambda $.
Then $s \in \{ 252 , 378 \} $ and
$\min (\Lambda ) \min (\Lambda ^* ) \in \{ 6 , \frac{14}{3} \} $.
\end{theorem}

The next two sections treat the two possibilities $s=252$ and $s=378$ separately.

\subsection{The case $s=252$.}

In this section we show

\begin{theorem}\label{s252}
There is no strongly perfect lattice in dimension $12$ with $s=252$.
\end{theorem}

Let $\Lambda $ be such a strongly perfect lattice with
 $s(\Lambda ) = s = 252 $ scaled such that $\min (\Lambda ) =2$ and
 put $\Gamma := \Lambda ^*$.
 Then $\min (\Gamma  ) $ is one of $3$ or $\frac{7}{3}$ and
 for all $\alpha \in \Gamma $ equation $(D4)$ yields that 
$
\sum _{x\in X }(x,\alpha )^4  = 18 (\alpha , \alpha )^2 
$
Hence $(\alpha ,\alpha) \in \frac{1}{3} \Z $. 
The equalities 
$\frac{1}{6} ((D13) - (D11))$, $\frac{1}{12} ((D4)-(D2)) $, and 
$(D22)$ 
yield that 
$$
\begin{array}{ll}
(\alpha ,\beta) (3 (\alpha ,\alpha ) -7) & \in \Z \\
\frac{1}{2} (\alpha ,\alpha) (3 (\alpha ,\alpha ) -7) & \in \Z  \\
6 (2 (\alpha , \beta )^2 + (\alpha , \alpha )(\beta , \beta )) & \in \Z 
\end{array}
\mbox{ for all } \alpha , \beta \in \Gamma .
$$
In particular, if $(\alpha ,\alpha )\in 2 \Z $, then 
$(\alpha , \beta ) \in \Z $ for all $\beta \in \Gamma ^* $.

Put 
$$\Delta := 
\{ \alpha \in \Gamma \mid (\alpha , \alpha ) \in 2\Z \} .$$
Then $\Delta $ is an even sublattice of index $\iota := [\Gamma : \Delta ]$
in $\Gamma $.
Moreover $\min (\Gamma ^*) = 2$ yields the bound
 $$\det (\Gamma ) \leq 16.17 ,\ \ \det(\Delta ) \leq 16.17 \iota ^2 .$$
 In particular $\min (\Delta )  = 4$.

\begin{lemma}\label{RR}
$\iota := |\Gamma / \Delta | $ is a divisor of  $12$.
If $\iota \neq 2, 6$ then the quadratic group  $(\Gamma / \Delta , -(,) ) $ is 
isometric to $R^*/R$ for the following root lattice $R$:
$$
\begin{array}{|l|c|c|c|}
\hline
\iota  & 3& 4 &  12 \\ 
\hline
R &    A_2 & D_4  & A_2 \perp D_4 \\
\hline
\end{array}
$$
If $\iota =2 $ or $\iota = 6 $ then 
$$\Gamma ^{(t)} := \{ \alpha \in \Gamma \mid (\alpha , \alpha ) \in \Z \} $$
is an odd integral sublattice of $\Gamma $ of index $\iota /2 $ 
such that $\Delta $ is the even sublattice of $\Gamma ^{(t)}$.
For $\iota = 6$ the quadratic group 
$\Gamma /\Gamma ^{(t)} $ is isometric to $(A_2^*/A_2 , - (,)) $.
\end{lemma}

\bew
By Lemma \ref{index4} and \ref{index3} the index of $\Delta $ in $\Gamma $
is a divisor of $12$.
Let us treat the primes 2 and 3 separately.
Assume that $3$ divides $\iota $.
Then there is $\alpha \in \Gamma $  with $(\alpha ,\alpha ) \in \frac{1}{3} \Z \setminus \Z $.
Equation $\frac{1}{12} ((D4)-(D2) ) $ yields that 
$3 (\alpha , \alpha ) - 7  \in 3 \Z $ hence $(\alpha , \alpha ) \in \frac{1}{3} + \Z $.
Hence the 3-Sylow subgroup of the quadratic group $\Gamma / \Delta $ 
is isometric to $( A_2^* / A_2 , - (,) )$.
For the 2-Sylow subgroup note that the norms in $\Gamma $ are 
2-integral. 
Moreover  $\alpha + \Delta  = \beta + \Delta \neq \Delta $ 
if and only if $(\alpha , \alpha ) \equiv (\beta , \beta ) \equiv 1 \pmod{2\Z _2} $ and $(\alpha ,\beta ) \in \Z _2$.
Hence if $|\Gamma / \Delta | = 4 $ or $12$, then 
the 2-Sylow subgroup of the quadratic group $\Gamma /\Delta $ is
isometric to $D_4^*/D_4$.
Also if $|\Gamma / \Delta | = 2 $ or $6$, then $\Gamma $ is 
already $2$-integral.
\eb

The strategy is to construct an even overlattice $\tilde{\Delta }$ of
$\Delta $, such that $\Delta = \Gamma \cap \tilde{\Delta } $ and 
hence 
$$\Gamma /\Delta \cong (\Gamma + \tilde{\Delta } )/ \tilde{\Delta } \cong 
- R^*/ R .$$
Denote the last isomorphism by $\varphi $.
Then the subdirect product 
$$M:= (\Gamma + \tilde{\Delta }) \subdir{} R^* 
= \{ (v_1,v_2) \in 
 (\Gamma + \tilde{\Delta }) \perp R^*  \mid v_2 = \varphi (v_1) \} $$
 is an even lattice of dimension 
 $12 + \dim (R) $ and determinant $\det(M) = \det (\tilde{\Delta }) /\det (R) $.

\begin{lemma}\label{n2z}
Let $\gamma \in \Delta _4$. Then 
$N_2(\gamma ) = \{ x_1,\ldots , x_5, y_1,\ldots, y_5 \}$ with 
$\gamma =x_i + y_i $, $( x_i,y_i ) = 0$ and $(x_i,y_j) = (x_i,x_j) = 1 $ for
all $i \neq j \in \{ 1,\ldots , 5 \}$.
The lattice generated by $N_2(\gamma )$ is isometric to the root lattice
$D_6$.
\end{lemma}

\bew
By Lemma \ref{linkomb} $|N_2(\gamma ) | = 10$.
Moreover if $x\in N_2(\gamma )$, then also $\gamma -x \in N_2(\gamma )$.
This gives us the partition 
$N_2(\gamma ) = \{ x_1,\ldots , x_5, y_1,\ldots, y_5 \}$ with  $x_i+y_i = \gamma $.
Taking scalar products with $x_1$ we get $$2 = (x_1,\gamma ) = (x_1,x_j) + (x_1,y_j).$$
Since $ (x_1,x) \leq 1$ for all $x\in N_2(\gamma ) $, $x\neq x_1$, we get
$(x_1,x_j) = (x_1,y_j) = 1$ for all $j\geq 2$ as claimed.
\eb

\begin{proposition}\label{indexvonDelta}
Fix some $\gamma \in \Delta _4$ and put 
$$\tilde{\Delta }:= \langle \Delta , N_2(\gamma  ) \rangle .$$
Then $[\tilde{\Delta } : \Delta ]  \geq 8$.
If $\iota = 3 $ or $\iota  = 12  $ then
 $[\tilde{\Delta } : \Delta ]  = 8$ 
and 
$\tilde{\Delta }/\Delta $ is not cyclic.
\end{proposition}

\bew
First we note that $\tilde{\Delta }$ is an even lattice,
$\Gamma \cap \tilde{\Delta } = \Delta $ and 
$\Gamma + \tilde{\Delta } \leq \tilde{\Delta }^*$.
Therefore $\iota $ divides the determinant of $\tilde{\Delta }$.

We keep the notation of Lemma \ref{n2z} and denote by 
$x\mapsto \overline{x} $ the natural projection  $\tilde{\Delta } \to \tilde{\Delta }/\Delta $.
Let $J:=\{ i\in \{ 1,\ldots , 5\} \mid 2\overline{x_i } \neq 0 \}$ and 
$I:=\{ 1,\ldots , 5\} \setminus J$.
Since $\min (\Delta ) = 4 > 2 $, 
we have $\overline{x_i} \neq 0 \neq \overline{x_i-x_j }$ for $i\neq j$. Hence
$\{ 0, \overline{x_1 } , \ldots , \overline{x_5 } \}$
are 6 distinct elements of $\tilde{\Delta }/\Delta $. Moreover, if 
$\overline{x_i} = -\overline{x_j }$ then 
$x_i +x_j -\gamma  = x_j -y_i \in \Delta $ is 
an element of norm $2$ for $j\neq i$.
Therefore
$\overline{x_i} = -\overline{x_j }$  if and only if $i=j \in I$. 
This implies that $ |\tilde{\Delta } / \Delta | \geq 6 + |J| $.

If $|\tilde{\Delta } / \Delta | = 6 $ then $J  = \emptyset $ and 
$\tilde{\Delta }/\Delta $ is an elementary abelian $2$-group.
Therefore $|\tilde{\Delta } /\Delta | \geq 8$.
If $|\tilde{\Delta }  / \Delta | = 7$ or $9$, then 
$\tilde{\Delta }/\Delta $ has no elements of order 2, hence $I=\emptyset $
and therefore $|\tilde{\Delta }/\Delta | \geq 11$.
Therefore 
$$|\tilde{\Delta}/\Delta | = 8 \mbox{ or } 
|\tilde{\Delta}/\Delta |  \geq 10 .$$

Now assume that $\iota \in \{ 3,  12 \}$ and let 
$R$ be the root lattice from Lemma \ref{RR}.
If $|\tilde{\Delta }/\Delta | \geq 10 $ then
$\det(\tilde{\Delta }) \leq \frac{\det (\Delta) }{100} \leq 23.3 \iota / 12 $.
Since $\det(\tilde{\Delta }) $ is divisible by $\iota $, we get 
$\det (\tilde{\Delta }) = \iota $ and 
$$ \tilde{\Delta }^* = \Delta + \Gamma  .$$
Taking the subdirect product 
$M:= (\Gamma + \tilde{\Delta }) \subdir{} R^* $ as described above,
constructs an even unimodular lattice in dimension 14 or  18,
which is a contradiction.
\eb

\begin{proposition}
$\iota = [ \Gamma : \Delta ] \neq 2$ or $6$.
\end{proposition}

\bew
Fix some $\gamma \in \Delta _4$ and put                                         $$\tilde{\Gamma }:= \langle \Gamma ^{(t)} , N_2(\gamma  ) \rangle .$$
As in the proof of Proposition \ref{indexvonDelta} we see that 
$[\tilde{\Gamma } : \Gamma ^{(t)} ] \geq 8$.

Assume that $\iota = 2$. Then $\Gamma  = \Gamma ^{(t)} $ is an integral lattice
of determinant  $\leq 16.17 $.
Hence $\det (\tilde{\Gamma }) \leq \frac{1}{8^2} 16.17  < 1 $ 
which is a contradiction.

If $\iota = 6$, then $\Gamma / \Gamma ^{(t)} \cong - A_2^*/ A_2 $.
Taking the subdirect product 
$ \Gamma \subdir{} A_2^* $ we obtain a 14-dimensional integral lattice of
determinant 
$\leq \frac{3 }{8^2} 16.17  < 1  $ which is a contradiction.
\eb

\begin{definition}
Let 
$$M:= (\Gamma + \tilde{\Delta }) \subdir{} R^* 
= \{ (v_1,v_2) \in 
 (\Gamma + \tilde{\Delta }) \perp R^*  \mid v_2 = \varphi (v_1) \} .$$
 \end{definition}

 \begin{remark}
  $M$
 is an even lattice of dimension 
 $12 + \dim (R) $ and determinant 
 $$\det(M) = \frac{\det (\tilde{\Delta }) }{\det (R)}  \leq
  \frac{\det (\Delta ) }{8^2 \det (R)} = \iota \frac{\det(\Gamma )}{64} 
  \leq 0.253 \iota .
  $$
  In particular $\iota \geq 4$ and hence $\iota \neq 3$.
  \\
 Moreover $R ^{\perp } =\{ m\in M \mid (m,r) = 0 \mbox{ for all } r\in R \}
 = \tilde{\Delta }$ contains a sublattice $R' = \langle N_2(\gamma ) \rangle
 \cong D_6$.
  \end{remark}

\begin{lemma}\label{component}
The lattice $R'$ is either a component of the root sublattice
$\langle \tilde{\Delta }_2 \rangle$ or $R' \subseteq E_7 \subseteq \langle 
\tilde{\Delta }_2 \rangle $
is contained in a component isometric to $E_7$.
\end{lemma}

\bew
Let $\tilde{R}$ be the component of $\langle \tilde{\Delta }_2 \rangle $
containing $R'$. 
We have 
$$\gamma = x_1+y_1 \in R' \subset \tilde{R} \subset \tilde{\Delta } \subset \Lambda  .$$
Since 
$N_2(\gamma )  \subset R'$,
there are no $ x\in \tilde{R}_2 \setminus R' $  satisfying $ (\gamma , x ) = 2$.
Going through all possible root lattices, we see that
only the two possibilities $\tilde{R} = R' \cong D_6$ or
$\tilde{R} \cong E_7$ arise.
\eb

We treat the two remaining cases $\iota = 4 $ and $\iota  = 12 $ separately.

\begin{lemma}
$\iota \neq 4$.
\end{lemma}

\bew
If $\iota = 4 $ then $M$ is a 16-dimensional even unimodular lattice, 
hence $M\cong E_8 \perp E_8$ or $M\cong D_{16}^+$ 
and $\tilde{\Delta } = D_4^{\perp }$ is the orthogonal complement of
a root system $D_4$ in $M$ and hence has root system
$E_8\perp D_4$ respectively $D_{12}$ contradicting Lemma \ref{component}.
\eb

%

\begin{lemma}
$\iota \neq 12$.
\end{lemma}

\bew
If $\iota =12$, then $M$ is an even 18-dimensional lattice of determinant
$\leq 3$.
Since there are no even unimodular 18-dimensional lattices  and also
no even lattices of determinant 2 and dimension 18, we have 
$\det (M) = 3$, $\det (\tilde{\Delta }) = 36$ and 
$\det(\Delta ) = 8^2\cdot 3 \cdot 12 = 2304 $.
There is one genus of 18-dimensional even lattices of determinant 3,
it contains 6 isometry-classes, representatives 
$M_1,\ldots , M_6$ of which have root sublattices 
$$E_7 \perp D_{10}, \  E_8\perp E_8 \perp A_2,  \ 
A_2 \perp D_{16}, \  A_{17}, \  E_6\perp E_6 \perp E_6,  \mbox{ and } D_7 \perp A_{11} .$$
By Lemma \ref{component}, the only possibility is $M = M_1$ with
root system $E_7\perp D_{10}$.
Since all reflections along norm 2 vectors are automorphisms of $M$,
there are up to isometries 3 embeddings of 
$R= A_2\perp D_4 $ into $M$:
$$ (1) A_2 \perp D_4 \subset D_{10} , \ \ 
(2) A_2 \subset D_{10}, \ D_4 \subset E_7 , \ \
(3) D_4 \subset D_{10} , \ A_2 \subset E_7 $$
(see for instance \cite[Table 4]{King}).
The root systems of the orthogonal complements of $R$ are
$$(1) A_3 \perp E_7 , \ \ (2) A _1^3 \perp D_7, \ \ (3) A_5 \perp D_6 .$$
By Lemma \ref{component} possibility (2) is excluded.
For the other two lattices we calculate all sublattices of 
index 8, such that the factor group is not cyclic to get 
a list of candidates for $\Delta $.
None of these lattices has minimum 4. 
\eb

This concludes the proof of Theorem \ref{s252}.

\subsection{The case $s=378$}

This section completes the proof of Theorem \ref{main12} by showing the following

\begin{theorem}\label{s378}
If $\Lambda $ is a strongly perfect lattice of dimension $12$ with 
Kissing number $2\cdot 378 $ then $\Lambda $ is similar to the
Coxeter Todd lattice $CT$.
\end{theorem}

Let $s=2\cdot 3^3 \cdot 7 = 378$ and rescale the strongly perfect
$12$-dimensional lattice $\Lambda $ such that $\min (\Lambda ) = \frac{4}{3}$.
Because of Remark \ref{CTmin}
we may and will assume that $\Lambda $ is generated by its minimal
vectors.
Let $\Gamma := \Lambda ^*$.
Then $\min (\Gamma ) = 4$ or $\min (\Gamma ) = \frac{7}{2}$
and 
for all  $\alpha, \beta \in \Gamma $
$$ 
\begin{array}{lrl}
(D2) &  2\cdot 3 \cdot 7 (\alpha  , \alpha  )  & \in \Z \\
(D4) &   2^2 \cdot 3 (\alpha  , \alpha  )^2  & \in \Z \\
12^{-1} ((D4)-(D2)) &  \frac{1}{2} (\alpha , \alpha ) (2(\alpha , \alpha ) - 7) & \in \Z  \\
6^{-1} ((D13)-(D11)) &  (\alpha , \beta ) (2(\alpha , \alpha ) - 7)
 & \in \Z  
 \end{array}
$$

From $(D2)$ we find that  $(\alpha ,\alpha ) \in \frac{1}{2} \Z $ 
for all $\alpha \in \Gamma $.
If $(\alpha ,\alpha ) \in \Z $, then 
$\frac{1}{6} ((D13)-(D11))$ implies that 
$(\alpha , \beta ) \in \Z $ for all $\beta \in \Gamma $.
In particular $(2\alpha , \beta ) = 2 (\alpha , \beta ) \in \Z $ 
for all $\alpha ,\beta \in \Gamma $.
Let $$\Gamma ^{(e)} := \{ \alpha \in \Gamma \mid (\alpha , \alpha ) \in \Z \} .$$ 
Then $\Gamma ^{(e)}$ is an integral sublattice of $\Gamma $ and
$\iota := [\Gamma : \Gamma ^{(e)} ] \leq 2$ since 
$(\alpha , \beta ) \in \frac{1}{2} \Z $ for all $\alpha ,\beta \in \Gamma $
and therefore $\Gamma ^{(e)} $ is the kernel of the
linear mapping $\Gamma \to \Z/2\Z, \alpha \mapsto 2(\alpha , \alpha )$.

Equality $\frac{1}{12} ((D4)-(D2)) $ yields that 
$(\alpha, \alpha )$ is even,
whenever $\alpha \in \Gamma ^{(e)}$, hence $\Gamma ^{(e)}$ is an even lattice.
Moreover the norms of the elements $\alpha \in \Gamma  \setminus \Gamma ^{(e)}$
satisfy $2 (\alpha , \alpha )  - 7 \in 4 \Z $, hence 
$$(\alpha , \alpha ) \in \frac{3}{2}  + 2 \Z \mbox{ for all }
\alpha \in \Gamma \setminus \Gamma ^{(e)} .$$
Moreover 
$$(\Gamma ,\Gamma ^{(e)}) \subseteq \Z , \ 
(\Gamma \setminus \Gamma ^{(e)} , \Gamma \setminus \Gamma ^{(e)} ) \in \frac{1}{2} +\Z $$
where the latter follows since $\alpha + \beta \in \Gamma ^{(e)}$ for all
$\alpha , \beta \in \Gamma \setminus \Gamma ^{(e)}$.
By Lemma \ref{min} we get that 
$$6.4 \geq \gamma _{12}^2 \geq \min (\Lambda ) \min (\Gamma )  
= \frac{4}{3} \min (\Gamma ) \geq \frac{14}{3}  $$
and therefore $\min (\Gamma ) \in \{ \frac{7}{2}, 4 \}$.
Let $d:=\det (\Gamma )$ and $d_0 := \det (\Gamma ^{(e)}) = \iota ^2 d$
for $\iota = [\Gamma : \Gamma ^{(e)} ]  \in \{ 1, 2 \}$. 
Then $$\min (\Lambda ) \sqrt[12]{d}  
= \frac{4}{3} \sqrt[12]{d_0/\iota } \leq \gamma _{12} \leq 2.522 $$
which yields $d_0\leq 2097 \iota^2 $.
On the other hand 
$\min (\Gamma ^{(e)}) \geq 4$ yields that $d_0 \geq 254$.
If $\min (\Gamma ^{(e)}) \geq 6$, then 
$d_0 \geq 32875 $ which is a contradiction.
Therefore there is some $\gamma \in \Gamma ^{(e)} $ with $(\gamma , \gamma ) = 4$.

Hence we have shown that

\begin{lemma}\label{minGamma0eq4}
$\Gamma ^{(e)} := \{ \alpha \in \Gamma \mid (\alpha , \alpha ) \in \Z \} $ 
is an even sublattice of $\Gamma $ of index $\iota \leq 2$ and minimum $4$.
We have $\det (\Gamma ) \leq 2097 $ and $\det (\Gamma ^{(e)} ) \leq 2097 \iota ^2  $.
\end{lemma}

The strategy of the proof of Theorem \ref{s378} 
is to construct an integral overlattice $M$ 
of $\Gamma ^{(e)}$  (respectively glue $\Gamma $ with $A_1^*$ to
obtain an even overlattice  of $\Gamma ^{(e)} \oplus A_1 $ in dimension
13 and then to find some overlattice $M$) of small determinant.
Then we go through all possibilites for $M$ and calculate $\Gamma ^{(e)}$
as a sublattice of $M$. 
To prove the existence of such an integral overlattice $M$, 
we find vectors of even norm and integral scalar products as 
linear combinations of the vectors in $X$ by analysing the 
possibilities for $N_2(\gamma )$ for $\gamma \in \Gamma $ satisfying the conditions
of Lemma \ref{linkomb}. 
The design properties of $X$ allow to obtain very precise information 
about the sets $N_2(\gamma )$, if $\gamma $ has norm 4 
(see Propositions \ref{N2gamma} and \ref{classes}).
Therefore we want to show that $\Gamma _4$ spans a space of dimension 
at least 8. This allows to construct an overlattice $M$ (see Corollary \ref{overlat}).
To use the theorem by Minkowski on the successive minima of $\Gamma $
we need to bound the number of vectors of norm $\frac{7}{2}$ in $\Gamma $.
Proposition \ref{gamma7/2}  shows that there are at most two such 
vectors and Proposition \ref{gamma11/2} shows that no vector in $\Gamma $ 
has norm $\frac{11}{2} $.

\begin{proposition}\label{gamma11/2}
There is no $\gamma \in \Gamma $ with $(\gamma ,\gamma ) = \frac{11}{2} $.
\end{proposition}

\bew
Let $\gamma \in \Gamma $ with $(\gamma ,\gamma ) = \frac{11}{2} $.
Then  by Lemma \ref{linkomb} $N_2(\gamma ) = \{ x_1 , \ldots , x_{11} \}$
with $\sum _{i=1}^{11} x_i = 4 \gamma $.
Since 
$(x_i,\sum _{j=1}^{11} x_j) = 4 (x_i,\gamma ) = 8$ and
$(x_i,x_j )\leq \frac{2}{3} $ for all $i\neq j$ we get that 
$(x_i,x_j ) = \frac{1}{10} (8-\frac{4}{3}) = \frac{2}{3} $ for all $i\neq j$ and 
$L:=\langle x_1,\ldots , x_{11} \rangle  \cong \sqrt{\frac{2}{3}} A_{11 } $.
We now enlarge $\Gamma ^{(e)}$ to an integral overlattice 
$$\Gamma ^{(e)} \subset \tilde{\Gamma } \subset \langle \Gamma ^{(e)} , N_{2}( \gamma ) \rangle $$
by joining preimages of a maximal isotropic subspace of 
$L/(L^*\cap L) \otimes \F _3 $.
We find such a subspace of dimension 5 which allows us to construct an
integral overlattice of $\Gamma ^{(e)}$ of index $3^5$ which contradicts the fact that 
$\det (\Gamma ^{(e)}) \leq 4 \cdot 2097 < (3^5)^2 $.
In detail let 
$$\begin{array}{ll}
y_1 : =  & x_3+x_4+x_5  \\
y_2 : =  & x_6+x_7+x_8  \\
y_3 : =  & x_9+x_{10}+x_{11}  \\
y_4 := & x_3 - x_4 + x_7-x_8+x_9-x_{11} \\
y_5 := & x_1 +x_7-x_8+x_{10} - x_{11} 
\end{array}
$$
Then the subspace $<y_1,\ldots , y_5 > \leq L$ is an integral 
sublattice of $L$ and
the linear functionals $x\mapsto (x,y_i) \in \frac{1}{3} L^*$
$(i=1,\ldots , 5)$ are linearly independent modulo $L^*$.
Therefore $$\tilde{\Gamma } := \langle \Gamma ^{(e)} , y_1,\ldots , y_5 \rangle $$
is an integral overlattice of $\Gamma ^{(e)}$ of index divisible by $3^5$.
On the other hand $$\det(\tilde{\Gamma } ) \leq \frac{1}{3^{10}} \det (\Gamma ^{(e)}) 
\leq \frac{4\cdot 2097}{3^{10}} \leq 0.15 < 1 $$
yields a contradiction.
\eb

The next aim is to investigate the vectors of norm 4 in $\Gamma $. 
From Lemma \ref{linkomb} and the equalities $(D2)$ and $(D4)$ we find

\begin{lemma}{\label{norm4}}
Let $\gamma \in \Gamma ^{(e)}$ with $(\gamma , \gamma ) =4$.
Then 
$N_2(\gamma ) = \{ x_1,x_2 \}$ with 
$(x_1,x_2) = \frac{2}{3} $ and $\gamma = x_1+x_2 $.
Moreover
$|N_1(\gamma )| = 160$ and $|N_0(\gamma ) | = 216 $.
\end{lemma}

\begin{proposition}\label{N2gamma}
Let $\gamma \in \Gamma ^{(e)}$ with $(\gamma , \gamma ) =4$ and fix
$x_1\in N_2(\gamma ) = \{ x_1,x_2\}$.
Then $$(x_1,X) \subset \{ 0,\pm \frac{1}{3} , \pm \frac{2}{3} , \pm \frac{4}{3} \}  $$
with 
$$
|N_0(x_1)| = 135 , \ 
|N_{1/3}(x_1)| = 80 , \ 
|N_{2/3}(x_1)| = 82 , \ 
|N_{4/3}(x_1)| = 1. 
$$
Choose $X$ such that $(x,\gamma ) \geq 0$ for all $x\in X$ and
$(x,x_1) \geq 0 $ for all $x\in N_0(\gamma )$ and define
$$M_{ij}(\gamma , x_1):= \{  x\in X \mid (x,\gamma ) =  i, (x,x_1) =  j/3 \}.$$
Then 
$M_{0,2}(\gamma , x_1) = \{ x_1 - x_2 \}$, $M_{2,2} (\gamma , x_1)= \{ x_2 \}$,
$M_{2,4} (\gamma , x_1) = \{ x_1 \}$ and the 
$m_{ij} := |M_{ij} (\gamma , x_1 ) |$ are given in the following table:
$$
\begin{array}{|c|c|c|c|c|c|} 
\hline
m_{ij}    & j=0 & j=1 & j=2 & j=4  
\\ \hline
i=0  & 135 & 80 & 1 & 0 
\\ \hline
i=1  & 0 & 80 & 80 & 0 
\\ \hline
i=2  & 0 & 0 & 1 & 1
\\ \hline
\end{array}
$$
\end{proposition}

\bew
Choose $X$ such that $(x,\gamma )\geq 0 $ for all $x\in X$.
Then 
$$X = X_0 \cup X_1 \cup X_2 $$
where $X_i = \{ x\in X \mid (x,\gamma ) = i \} $ and
$X_2 = \{ x_1,x_2 \}$ by Lemma \ref{norm4}.
We also choose the elements in $X_0$, such that $(x,x_1) \geq 0 $ for all $x\in X_0$.
Then the equation 
$(D11)$ with $\alpha _1 = \gamma $ yields that for all $\alpha \in \R^{12}$
$$\sum _{ x\in X_1 } (x,\alpha ) + 2(x_1,\alpha ) + 2 (x_2,\alpha ) =
2\cdot 3 \cdot 7 (\gamma , \alpha )  .$$
Since $\gamma = x_1 + x_2 $ this gives
$$\sum _{ x\in X_1 } (x,\alpha ) 
= 2^ 3 \cdot 5 (\gamma , \alpha )  .$$
Similarly equality $(D22)$ yields
$$\sum _{x\in X_1} (x,\alpha )^2 = 2^3 (\gamma ,\alpha )^2
-2^2 (x_1,\alpha )^2 - 2^2 (x_2,\alpha )^2 + 2^4 (\alpha , \alpha ) $$
and 
$(D13)$ gives 
$$\sum _{x\in X_1} (x,\alpha )^3 = 2^2 3 (\gamma ,\alpha ) (\alpha , \alpha )
-2 (x_1,\alpha )^3 - 2 (x_2,\alpha )^3  $$

Let $x\in X_1$. Then $1=(x,\gamma ) = (x,x_1 ) + (x,x_2)$.
Since $(x,x_i)\leq \frac{2}{3}$, this implies that 
$(x,x_1) \geq \frac{1}{3}$.
For $i\in \R $ let $m_i := | \{ x\in X_1 \mid (x,x_1) = i \} |$.
Then the equalities above yield
$$
\begin{array}{ll}
\sum _{i} m_i &  = |X_1| = 2^5\cdot 5  \\
\sum _i i m_i  & =  2^4 \cdot 5 \\
\sum _i i^2 m_i  & =  2^4 \cdot 5^2 \cdot 3^{-2} \\
\sum _i i^3 m_i  & =  2^4 \cdot 5 \cdot 3^{-1} 
\end{array}
$$
Hence 
$$\sum _i (i-\frac{1}{3})(\frac{2}{3} -i)  m_i 
= \sum _i (-i^2+i-\frac{2}{9} ) m_i = 0.
$$
Since $\frac{1}{3} \leq i \leq \frac{2}{3} $ this yields that 
$(x,x_1) \in \{ \frac{1}{3}, \frac{2}{3} \}$ for all $x\in X_1$.
Moreover we find
$$m_{1,1} = 80, \ \ m_{1,2} = 80 .$$

Now equalities $(D2)$ and $(D4) $ yield the following equations 
for $n_i := \{ x\in X_0 \mid (x_1,x) = i \}$:
$$\begin{array}{ll}
\sum _{i} n_i & = |X_0 | = 2^3 3^3 \\
\sum  _{i} i^2 n_i  & = 2^2\cdot 7\cdot 3^{-1} \\
\sum _{i} i^4 n_i & = 2^5 \cdot 3^{-3} 
\end{array}
$$
By our assumption $(x,x_1) \geq 0 $ for all $x\in X_0$.
Moreover for $x\in X_0$ we have $(x,\gamma ) = (x,x_1+x_2) = 0$ hence 
$(x,x_2) = - (x,x_1)$.
Since $(x_1,x_2) = 2/3$ we get 
$(x_1-x_2,x_1-x_2) = 4/3$ and hence
$( x, x_1-x_2) = 2(x,x_1) \leq 2/3 $ for all $x\in X_0 $ with $x\neq x_1-x_2$.
Therefore $n_{2/3} = 1$ and $n_i \neq 0$ only for $i\in [0,\frac{1}{3} ] \cup \{ 2/3 \} $
and we get 
$$\begin{array}{ll}
\sum _{i \in [0,\frac{1}{3} ]} n_i & = |X_0 | - 1 = 2^3 3^3 -1 \\
\sum  _{i \in [0,\frac{1}{3} ]} i^2 n_i  & = 2^4 \cdot 5\cdot 3^{-2} \\
\sum _{i \in [0,\frac{1}{3} ]} i^4 n_i & = 2^4 \cdot 5\cdot 3^{-4} 
\end{array}
$$
from which we get
$$
\sum _{i \in [0,\frac{1}{3} ]} i^2 (i^2-\frac{1}{9}) n_i  =  0 .
$$
Therefore $n_i = 0 $ for $i \not\in \{ 0,\pm \frac{1}{3} , \pm \frac{2}{3} \} $
and we may use the 3 equalities above to calculate 
$$n_0 = 135, n_{1/3} = 80, n_{2/3} = 1 $$
from which the proposition follows.
\eb

\begin{corollary}
In the situation of Proposition \ref{N2gamma} we have 
$3x_1 \in \Gamma $ and $3x_2 \in \Gamma $.
Hence also $\gamma ' := x_1 -2x_2 = \gamma - 3 x_2 \in \Gamma $
with $(\gamma ' , \gamma ' ) = 4$ and $(\gamma , \gamma ' ) = -2 $.
\end{corollary}

This allows us to define an equivalence relation on the set of norm 4 vectors
$C :=  \Gamma _4$ in $\Gamma $.
Note that $C$ consists of the minimal vectors in the integral sublattice $\Gamma ^{(e)}$
of $\Gamma $.
In particular $C$ is not empty by Lemma \ref{minGamma0eq4} and 
$|(\gamma _1,\gamma _2) | \leq 2 $ for distinct elements $\gamma _i \in C $.
We call $\gamma _1,\gamma _2 \in C $ {\em equivalent}, if $\gamma _1-\gamma _2 \in 3 \Lambda $.
Let ${\cal K} $ denote the set of equivalence classes.
Then  ${\cal K}$ forms a root system   over $\Z [\zeta _3 ]$.
More precisely we have:

\begin{proposition}{\label{classes}}
\begin{itemize}
\item[(i)] For all $K\in {\cal K}$ we have $|K| = 3$.
\item[(ii)] If $K=\{ u,v,w \}  $ then 
$u +v + w = 0 $ and $(u ,v) = (u,w) = (v,w) = -2 $.
Moreover $N_2(u) = \{  \frac{1}{3} (u-v) , \frac{1}{3} (u-w) = \frac{1}{3} (2u+v) \} $.
\item[(iii)] If $K_1 \neq \pm K_2 \in {\cal K} $ then either 
$(k_1,k_2) = 0 $ for all $k_i \in K_i$ (in this case the classes are called {\em orthogonal},
$(K_1,K_2) := 0 $) or 
there is $\epsilon \in \{\pm 1 \}$ and
a mapping $\varphi : K_1 \to K_2$ with $(k_1 , \varphi (k_1 )) = 2\epsilon  $ and
$(k_1,k_2) = -\epsilon $ for all $k_1 \in K_1, k_2\in K_2\setminus \{ \varphi(k_1) \} $.
In this situation we will say that $(K_1,K_2) := \epsilon $.
\item[(iv)] If $(K_1,K_2) = -1$, then $K_3 := K_1 + K_2 := \{ k_1 + \varphi (k_1)\mid 
k_1\in K_1 \} \in {\cal K}$.
\end{itemize}
\end{proposition}

\bew
(i +ii) Let $u \in K\in {\cal K}$. Then $N_2(u ) =  \{ x _{u } , y_{u } \}$ and
$\{ u, u - 3 x_{u }, u - 3 y_{u }  \} \subset K$.
On the other hand let $u \neq v \in K$. Then $x:=\frac{1}{3} (u-v) $
is a non-zero vector in $\Lambda $ and hence has square length $\geq \frac{4}{3}$.
This implies that $(u,v )  = -2$ and $x\in N_2(\gamma )$.
\\
(iii)
Since the differences of the elements in $K_i$ lie in $3 \Lambda  = 3 \Gamma ^*$, the
scalar products $(k_1, k_2) \equiv (k_1' , k_2' ) \pmod{3}$ are congruent modulo 3 for
all $k_i,k_i' \in K_i , i=1,2 $.
Now  these scalar products are integers  of absolute value $\leq 2$  with 
$ \sum _{k_2\in K_2} (k_1,k_2) = 0 $ for all $k_1\in K_1$ which only leaves the
possibilities described in the proposition.
\\
(iv)
Clearly $K_3$ is again an equivalence class of norm 4 vectors in $\Gamma $.
\eb

A closer analysis of the proof of Proposition \ref{N2gamma} allows to 
define a normal subgroup of the automorphism group of $\Lambda $.
For $\Lambda = CT$ this is the representation of a 
subgroup of index 2 in $\Aut(CT)$ as complex reflection group.

\begin{proposition}
For all $K\in {\cal K}$  define an orthogonal mapping 
$s_K$ by $(s_K)_{|L} = - \id _{|L }$ and
$(s_K)_{|L^{\perp }} = \id _{|L^{\perp }}$, where 
$L = \langle K \rangle _{\R }$ is the vector space generated by $K$.
If $\Lambda $ is generated by its minimal
vectors, then $s_K  (\Lambda ) = \Lambda $ hence $s_K \in \Aut (\Lambda )$.
Moreover 
\begin{itemize}
\item[(i)] $s_K^2 = \id $.
\item[(ii)] $(s_{K_1} s_{K_2} ) ^2 = \id $ for all $K_1,K_2 \in {\cal K}$, $(K_1,K_2) = 0$.
\item[(iii)] $(s_{K_1} s_{K_2} ) ^3 = \id $ for all $K_1,K_2 \in {\cal K}$, $(K_1,K_2) \neq 0$.
\item[(iv)] The subgroup $\langle s_K \mid K\in {\cal K} \rangle 
\leq \Aut (\Lambda )$ is a normal subgroup of the automorphism group of 
$\Lambda $.
\end{itemize}
\end{proposition}

\bew
Let $Y:=X\cup -X = \Lambda _{\min} $. We only need to show that
$s_K(Y) = Y$. The remaining properties of the $s_K$ follow from direct 
calculations. In particular for all $g\in \Aut (\Lambda )$
the conjugate $s_{K}^g = s _{g(K)}$.
Since $\Aut(\Lambda ) = \Aut(\Gamma )$ permutes the elements of 
$K$, the $s_K$ generate a normal subgroup
of $\Aut (\Lambda )$.
Let $\gamma \in K \in {\cal K}$ and let $N_2(\gamma ) = \{ x_1,x_2 \}$.
Then $$\langle K \rangle_{\R } = \langle x_1, x_2 \rangle _{\R } = \langle \gamma , x_1 \rangle _{\R }.$$
First let $y\in Y$ with $(y,\gamma ) = 0$.
If $(y,x_1) = 0$, then $y\in K^{\perp}$ and $s_K(y) = y \in  Y$.
If $(y,x_1) = \frac{1}{3}$, then  $(y,x_2) = -\frac{1}{3} $, 
$y-\frac{1}{2} (x_1-x_2) \in K^{\perp }$,
and $s_K(y) = y  - (x_1-x_2) \in Y $.
Similarly, if $(y,x_1) = -\frac{1}{3}$, then $s_K(y) = y +(x_1-x_2) \in Y$.
If $(y,x_1) = \frac{2}{3}$, then $y=x_1-x_2 \in K$ by Proposition \ref{N2gamma}.
Therefore $s_K(y) = -y \in Y$.
\\
If $(y,\gamma ) = \pm 2$, then $y\in \pm N_2 (\gamma ) \subset \langle K \rangle$ and hence
$s_K(y) = -y  \in Y$.
\\
It remains to consider the case that $(y,\gamma ) = \pm 1$.
Without loss of generality let $(y,\gamma  ) = 1$. 
Then by Proposition \ref{N2gamma} $(y,x_1)$ is one of $\frac{1}{3}$ or $\frac{2}{3}$.
In the first case, the projection of $y$ onto $\langle K \rangle _{\R }$ 
is $\frac{1}{2} x_2$ and hence $s_K(y) = y - x_2 \in Y$.
In the second case, the projection of $y$ onto $\langle K \rangle _{\R }$ 
is $\frac{1}{2} x_1$ and hence $s_K(y) = y-x_1 \in Y$.
\eb

\begin{proposition}\label{a2a2}
Let ${\cal K}$ denote the set of equivalence classes introduced in 
Proposition \ref{classes}.
If there are two classes $K_1, K_2 \in {\cal K}$ with
$(K_1,K_2)= -1$ then $\Gamma \cong CT$ is the Coxeter-Todd lattice.
\end{proposition}

\bew
Let $K_i:=\{ u_i,v_i,-(u_i+v_i) \} $ 
($i=1,2$) such that the Gram matrix of
$(u_1,v_1,u_2,v_2) $ is 
$$A:= \left( \begin{array}{rrrr} 4 & -2 & -2 & 1 \\
-2 & 4 & 1 & -2 \\
-2 & 1 & 4 & -2 \\
1 & -2 & -2 & 4 \end{array} \right) .$$
Then $N_2 (u_i) = \{ x_i : = \frac{1}{3} (u_i-v_i), y_i:= \frac{1}{3} (2u_i + v_i ) \} $  ($i=1,2 $)
and $t:= u_1 - v_2 $ is a vector of norm 6 in $\Gamma $.
We now want to investigate $N_2(t)$:
The elements $y_1, x_2-y_2,$ and $ x_1+x_2$ have scalar product 2 with $t$
and satisfy $y_1+x_2-y_2+x_1+x_2 = t$.
With Lemma \ref{linkomb} we find 
$N_2(t) = \{ z_1,\ldots , z_{12}, z_{13}:=y_1,z_{14}:=x_2-y_2,z_{15}:=x_1+x_2 \}$ with $\sum _{i=1}^{15} z_i = 5 t $.
Since 
$t=u_1-v_2$ we get 
$$ 2 = (z_i,t) = (z_i,u_1) - (z_i,v_2) = 1 - (-1) \mbox{ for all } i \in\{ 1,\ldots , 12 \} $$
Hence $\{z_1,\ldots , z_{12} \} \subseteq N_1(u_1) \cap N_{-1} (v_2 ) $.
Since $|(z_i,z_j) | \leq \frac{2}{3} $ for all $i\neq j$ and 
$(z_{13},z_i) + (z_{14}, z_i) + (z_{15} , z_i) = (t,z_i) = 2 $ for all $i \leq 12$, we find that 
$(z_j,z_i) = \frac{2}{3} $ for all $i\in \{ 1,\ldots , 12 \} $, $j\in \{ 13,14,15 \}$.
For $i=1,\ldots , 12$ let 
$\overline{z}_i := z_i - \frac{1}{3} t $. 
Then 
$$(\overline{z}_i , u_1) = (\overline{z}_i , v_2 ) =  (\overline{z}_i, z_j) = 0 
\mbox{ for all } 1\leq i \leq 12, j = 13,14,15 .$$
Hence $\{ z_1,\ldots , z_{12} \} \in \langle u_1, v_1,u_2,v_2 \rangle ^{\perp }$ lie in an 8-dimensional space.
Moreover $$(\overline{z} _i ,\overline{z} _j )  = (z_i,z_j) -\frac{2}{3} 
\left\{ \begin{array}{ll} =\frac{2}{3} & i = j 
\\ 
\leq 0 & i \neq j \end{array} \right. .$$
By Lemma \ref{neg} there is a partition 
$\{ 1,\ldots , 12\} = I_1\cup \ldots \cup I_k$ 
into disjoint sets $I_j$ such that 
the minimal relations among the $\overline{z}_i$ are of the
form $\sum _{i\in I_j} \overline{z}_i = 0$ for all $j=1,\ldots , k$.
Clearly $|I_j| > 1$.
If $|I_j| = 2 $, then $\frac{2}{3} t = z_i + z_l $ for some $i\neq l$ and hence 
$2t = 3(z_i+z_l) \in 3\Lambda $.
Since $t\in \Lambda $ this implies that $t=3t-2t \in 3\Lambda $, hence
$\frac{t}{3} \in \Lambda $ is a vector of norm $\frac{2}{3}$ contradicting the fact that
$\min (\Lambda ) = \frac{4}{3} $.
Therefore $|I_j| \geq 3$ for all $j=1,\ldots , k$.
Since $k\geq 4$ and the $I_j$ are disjoint, this implies that 
$k=4$ and $|I_j| = 3 $ for all $j$.
Rearranging the $z_i$, we may assume that $I_j = \{ 2j-1 ,2j ,8+j \} $ for $j=1,\ldots , 4$.
Put $$\begin{array}{ll}
a_j  & := (\overline{z}_{2j-1} , \overline{z}_{2j} ) \\
b_j  & := (\overline{z}_{2j-1} , \overline{z}_{8+j} ) \\
c_j  & := (\overline{z}_{2j} , \overline{z}_{8+j} ) .  \end{array} $$
Then $a_j+b_j=a_j+c_j = b_j+c_j = \frac{-2}{3} $ which implies that $a_j=b_j=c_j = \frac{-1}{3} $
for all $j$.
Hence $\Lambda $ contains a sublattice  
$L:= \langle u_1, t, x_2-y_2,x_1+x_2, z_1,\ldots , z_{8}  \rangle $ where the
Gram matrix of $L$ is
$$F:=\frac{1}{3} \left( \begin{array}{cccccccccccc}
 12 & 9 & 0 & 3 & 3 & 3 & 3 & 3&  3 & 3&  3&  3 \\
 9 & 18 & 6 & 6 & 6 & 6 & 6 & 6&  6 & 6&  6&  6 \\
 0 & 6 & 4 & 1 & 2 & 2 & 2 & 2&  2 & 2&  2&  2 \\
 3 & 6 & 1 & 4 & 2 & 2 & 2 & 2&  2 & 2&  2&  2 \\
 3 & 6 & 2 & 2 & 4 & 1 & 2 & 2&  2 & 2&  2&  2 \\
 3 & 6 & 2 & 2 & 1 & 4 & 2 & 2&  2 & 2&  2&  2 \\
 3 & 6 & 2 & 2 & 2 & 2 & 4 & 1&  2 & 2&  2&  2 \\
 3 & 6 & 2 & 2 & 2 & 2 & 1 & 4&  2 & 2&  2&  2 \\
 3 & 6 & 2 & 2 & 2 & 2 & 2 & 2&  4 & 1&  2&  2 \\
 3 & 6 & 2 & 2 & 2 & 2 & 2 & 2&  1 & 4&  2&  2 \\
 3 & 6 & 2 & 2 & 2 & 2 & 2 & 2&  2 & 2&  4 &  1  \\
 3 & 6 & 2 & 2 & 2 & 2 & 2 & 2&  2 & 2&  1 &  4
\end{array} \right)
$$
Therefore $\Gamma $ is a sublattice of the dual lattice $M:=L^*$ with
$\min (\Gamma ) = 4$.
In particular $\Gamma  = \Gamma ^{(e)}$, since $M$ is already an even lattice.
The determinant of $M$ is $\det (M) = 81$, hence 
$$[M:\Gamma ] \leq \sqrt{\frac{2097}{81}} \leq 5.1 .$$

Computations with MAGMA yield that 
$\Aut (M)$ has 6 orbits on the sublattices $M^{(2)}$ of index 2 in $M$, none of 
 which satisfies $\min (M^{(2)*}) \geq \frac{4}{3} $.
$\Aut (M)$ has $32$ orbits on the sublattices $M^{(3)}$ of index 3 in $M$.
only one
of which satisfies $\min (M^{(3)*}) \geq \frac{4}{3} $.
This lattice $M^{(3)}$ is isometric to the Coxeter Todd lattice $CT$.

$\Aut (M)$ has $253$ orbits on the sublattices  $M^{(5)}$ of index 5 in $M$.
For none of these sublattices $M^{(5)}$ the dual $M^{(5)*}$ 
  has minimum $\geq \frac{4}{3}$.
\eb

The strategy is now to give a lower bound on the rank of the sublattice of 
$\Gamma $ 
spanned by the norm $4$-vectors using Minkowski's theorem on the successive
minima of lattices.
To this aim, we want to bound the number of
vectors of norm $\frac{7}{2} $ in $\Gamma $.

Let  $A:=\{ \alpha \in \Gamma \mid (\alpha ,\alpha ) = \frac{7}{2} \} = \Gamma _{7/2} $ 
and $C:=\Gamma _4$ as above.

\begin{lemma}{\label{xxx1}}
Let  $\gamma \in C$ and denote $N_2(\gamma ):= \{ x_{\gamma }, y_{\gamma }\}$.
\begin{itemize}
\item[(i)]
For all $\alpha \in A$ we have 
$(\alpha , \gamma ) = (\alpha , x_{\gamma } ) = (\alpha , y_{\gamma }) = 0 $.
\item[(ii)]
For $\alpha _1 \neq \pm \alpha _2 \in A $ we have $(\alpha _1,\alpha _2) = \pm \frac{1}{2} $.
\end{itemize}
\end{lemma}

\bew
(i)
Since the scalar products  $(\Gamma , \Gamma ^{(e)})$ are integral, we have the 
possibilities 
$(\alpha , \gamma ) \in \{ 0,\pm 1, \pm 2 \}$.
If $(\alpha , \gamma ) = 1$, then $\alpha - \gamma \in \Gamma _{11/2} $ is a vector of 
norm $\frac{11}{2} $ in $\Gamma $ which is impossible by Proposition
\ref{gamma11/2}.
If $(\alpha , \gamma ) = 2$, then $(\alpha , x_{\gamma }) = (\alpha , y_{\gamma }) = 1$
(since both scalar products are $\leq 1$) and $\gamma -3 x_{\gamma } \in C$ has
scalar product $-1$ with $\alpha $ contradicting Proposition
\ref{gamma11/2}.
Therefore $(\alpha ,\gamma ) = 0$ for all $\alpha \in A$, $\gamma \in C$.
Since $\gamma = x_{\gamma } + y _{\gamma }$ either both 
 $(\alpha , x_{\gamma }) = (\alpha , y_{\gamma }) = 0$ or one of them is $1$.
 As above the latter allows to construct a vector of norm $\frac{11}{2} $ in $\Gamma $.
 \\
 (ii)
 The possible scalar products are 
 $(\alpha _1,\alpha _2 ) \in \{ \pm \frac{1}{2} , \pm \frac{3}{2}  \} $.
 If $(\alpha _1,\alpha _2) = \frac{3}{2} $, then $\gamma := \alpha _1 - \alpha _2 \in C$
 satisfies $(\gamma , \alpha _1) = 2$ contradicting (i).
\eb

\begin{proposition}\label{gamma7/2}
$\frac{1}{2} |\Gamma _{7/2} | \leq 1$.
\end{proposition}

\bew
Assume that  there are $\alpha _1 \neq \pm \alpha _2 \in \Gamma _{7/2}$.
Then by Lemma  \ref{xxx1} we have $(\alpha _1 , \alpha _2 ) = \pm \frac{1}{2} $.
Assume that  $(\alpha _1,\alpha _2) = - \frac{1}{2} $ and put
$t:= \alpha _1 + \alpha _2 \in \Gamma _6$.
With Lemma \ref{linkomb} we find 
$N_2(t) = \{ y_1,\ldots , y_{15} \}$ with $\sum _{i=1}^{15} y_i = 5 t$.
For $j=1,2$ we calculate 
$ 15 = 5(\alpha _j , t) = \sum _{i=1}^{15} (\alpha _j ,y_i )$ hence $(\alpha _j ,y_i ) = 1$
for all $i$.
For $i=1,\ldots , 15$ let 
$\overline{y}_i := y_i - \frac{1}{3} t $. Then $(\alpha _j ,\overline{y}_i ) = 0 $ for all $i\in \{ 1,\ldots, 15\} ,j = 1,2$ and
the $\overline{y}_i$ lie
in the 10-dimensional space $\langle \alpha _1,\alpha _2 \rangle ^{\perp }$.
Moreover $(\overline{y} _i ,\overline{y} _j )  = (y_i,y_j) -\frac{2}{3} 
\left\{ \begin{array}{ll} =\frac{2}{3} & i = j 
\\ 
\leq 0 & i \neq j \end{array} \right. $.
By Lemma \ref{neg} there is a partition 
$\{ 1,\ldots , 15\} = I_1\cup \ldots \cup I_k$ 
into disjoint sets $I_j$ such that 
the minimal relations among the $\overline{y}_i$ are of the
form $\sum _{i\in I_j} \overline{y}_i = 0$ for all $j=1,\ldots , k$.
As in the proof of Proposition \ref{a2a2} we find 
$k=5$ and $|I_j| = 3 $ for all $j$.
Rearranging the $y_i$, we may assume that $I_j = \{ 2j-1 ,2j ,10+j \} $ for $j=1,\ldots , 5$.
Put $$\begin{array}{ll}
a_j  & := (\overline{y}_{2j-1} , \overline{y}_{2j} ) \\
b_j  & := (\overline{y}_{2j-1} , \overline{y}_{10+j} ) \\
c_j  & := (\overline{y}_{2j} , \overline{y}_{10+j} ) .  \end{array} $$
Then $a_j+b_j=a_j+c_j = b_j+c_j = \frac{-2}{3} $ which implies that $a_j=b_j=c_j = \frac{-1}{3} $
for all $j$.
Hence $\Lambda $ is an overlattice of the lattice
$L:= \langle y_1,\ldots , y_{10} , t , 2 \alpha _1 \rangle $ where the
Gram matrix of $L$ is
$$F:=\frac{1}{3} \left( \begin{array}{cccccccccccc}
4 & 1 & 2 & 2 & 2 & 2 & 2 & 2&  2 & 2&  6&  6 \\
1 & 4 & 2 & 2 & 2 & 2 & 2 & 2&  2 & 2&  6&  6 \\
2 & 2 & 4 & 1 & 2 & 2 & 2 & 2&  2 & 2&  6&  6 \\
2 & 2 & 1 & 4 & 2 & 2 & 2 & 2&  2 & 2&  6&  6 \\
2 & 2 & 2 & 2 & 4 & 1 & 2 & 2&  2 & 2&  6&  6 \\
2 & 2 & 2 & 2 & 1 & 4 & 2 & 2&  2 & 2&  6&  6 \\
2 & 2 & 2 & 2 & 2 & 2 & 4 & 1&  2 & 2&  6&  6 \\
2 & 2 & 2 & 2 & 2 & 2 & 1 & 4&  2 & 2&  6&  6 \\
2 & 2 & 2 & 2 & 2 & 2 & 2 & 2&  4 & 1&  6&  6 \\
2 & 2 & 2 & 2 & 2 & 2 & 2 & 2&  1 & 4&  6&  6 \\
6 & 6 & 6 & 6 & 6 & 6 & 6 & 6&  6 & 6& 18 & 18  \\
6 & 6 & 6 & 6 & 6 & 6 & 6 & 6&  6 & 6& 18 & 42
\end{array} \right)
$$
satifying $\min (\Lambda ) = \frac{4}{3}$ and $\min (\Lambda ^*) = \frac{7}{2}$.
In particular $\Gamma ^{(e)}$ is an even integral sublattice of $L^*$ of minimum 4
and hence contained in the unique maximal integral sublattice $M_0$ of $L^*$.
If $B^*$ denotes the dual basis of the basis of $L$ above, then
$M_0$ is generated
by $(b_1^*,\ldots , b_{10}^*, 2 b_{11}^* + 2 b_{12}^* , b_{11}^*+3b_{12}^* )$
and has index 4 in $L^*$ and determinant $3^4$.
$M_0$ is an overlattice of index 3 of the root lattice $A_2^6$.
The lattice $\Gamma $ is  a sublattice of $M:=M_0 + \Z \alpha _1 $ of
index $\leq \sqrt{4\cdot 2097 / 81 }  \leq 10.2 $.

Computations with MAGMA yield that 
$\Aut (M)$ has 23 orbits on the sublattices $M^{(2)}$ of index 2 in $M$, only 
4 of which satisfy $\min (M^{(2)*}) \geq \frac{4}{3} $.
These sublattices have minimum $2,2,2,$ and $3/2$, hence are no
candidates for $\Gamma $.
All sublattices $M^{(4)}$ of these 4 sublattices $M^{(2)}$ of index 2 satisfy 
$\min (M^{(4)*})  < \frac{4}{3} $.
Hence $[M:\Gamma ] = [M_0 : \Gamma _0] $ is not a power of $2$.

Assume that $3$ divides $[M:\Gamma ]$.
$\Aut (M)$ has $109$ orbits on the sublattices $M^{(3)}$ of index 3 in $M$.
only one
of which satisfies $\min (M^{(3)*}) \geq \frac{4}{3} $.
Since $\min (M^{(3)}) = \frac{3}{2}$, the lattice $\Gamma $ is a proper sublattice of $M^{(3)}$.
There are no sublattices $M^{(9)}$ of index 3 of $M^{(3)}$ such that 
$M^{(9)*} $ has minimum $\geq \frac{4}{3}$.
A unique sublattice $M^{(6)}$ of index 2 in $M^{(3)}$ satisfies 
$M^{(6)*} = \frac{4}{3}$. This lattice is 
Coxeter-Todd lattice and has $\min (M_9) = 4 > \frac{7}{2}$.

Assume that $5$ divides $[M:\Gamma ]$.
$\Aut (M)$ has $1771$ orbits on the sublattices  $M^{(5)}$ of index 5 in $M$.
For none of these sublattices $M^{(5)}$ the dual $M^{(5)^*}$ 
has minimum $\geq \frac{4}{3}$.

It remains the case that $[M:\Gamma ] = 7$.
Here the orbit computations are too big to be performed with MAGMA.
We therefore calculate all sublattices of $M$ of index 7 by going through
the 11-dimensional subspaces of $M/7M \cong \F_7^{12}$.
These are parametrised by matrices in $\F_7^{11\times 12} $ that are in
Hermite normal form.
Building these matrices row by row, we only continue if the sublattice
generated by $7M$ and the first rows has minimum $\geq \frac{7}{2}$.
There are 31104 sublattices $M^{(7)}$ 
of $M$ of index 7 that have minimum $\geq \frac{7}{2}$ but
none of these lattices satisfies $\min (M^{(7)*}) \geq \frac{4}{3}$.
\eb

\begin{lemma}\label{mink7}
The rank of the sublattice of $\Gamma $ generated by the norm $4$ vectors
in $\Gamma $ is at least $7$.
\end{lemma}

\bew
We use Minkowski's theorem on the successive minima 
$$m_i:= \min \{ \lambda \in \R \mid L \mbox{ has } i \mbox{ linearly independent
vectors of norm } \leq \lambda \}$$
of an $n$-dimensional lattice $L$ stating that
$$m_1\cdot m_2 \cdot \ldots \cdot m_r \leq \gamma _n^r \det (L) ^{r/n} \ \mbox{ for all } 1\leq r\leq n  $$
(see for instance \cite[Th\'eor\`eme II.6.8]{Martinet}).
\\
Assume first that 
 $\Gamma $ does not contain vectors of norm $\frac{7}{2}$.
By Proposition \ref{gamma11/2}
there are no vectors of norm $\frac{11}{2}$ in $\Gamma $.
Since
$$4^6\cdot 6^6  > 2.522^{12} \cdot 2097 \geq \gamma _{12}^{12} \det (\Gamma ) $$
we get that the rank of the sublattice of $\Gamma $ spanned by $\Gamma _4$ is
at least 7. 
\\
Now assume that $\Gamma _{7/2} \neq \emptyset $.
Then by Proposition \ref{gamma7/2} the set $\Gamma _{7/2} = \{ \alpha , - \alpha \}$
has only 2 elements.
Let $\pi $ denote the orthogonal projection onto $\langle \alpha \rangle ^{\perp }$ and let $M:=\pi (N_1(\alpha ))$.
Then  by Lemma \ref{2design}, the set $\pi (M) \disj -\pi (M) $ is a 2-design
in $\langle \alpha \rangle ^{\perp }$.
In particular  $\langle \pi (M) \rangle = \langle \alpha \rangle ^{\perp }$
and hence $M$ spans a subspace of dimension $\geq 11$ of $\R ^{12}$.
\\
Choose $x\in M$ and let $G(x) := \{ \gamma \in \Gamma \mid (\gamma ,x ) = 0 \}$.
Then $G(x) $ is an 11-dimensional  lattice of determinant 
$$\det (G(x)) = \frac{4}{3} \det (\Gamma ) \leq 2796 .$$
Since $(\alpha , x) = \pm 1$, the minimum of $G(x) $ is $\geq 4$ and
$G(x) $ contains no vectors of norm $\frac{11}{2}$ by Proposition \ref{gamma11/2}.
Now  $\gamma_{11} \leq 2.39 $ by \cite{Elkies} and 
$$4^5\cdot 6^6  > 2.39^{11} \cdot 2796  $$
which implies that the rank of the sublattice spanned by the 
norm 4-vectors in $G(x)$  is at least 6.
This holds for any $x\in M$. Since $M$ spans a space of dimension $\geq 11$,
not all norm 4-vectors of $\Gamma $ can be orthogonal to all $x\in M$.
Therefore there is some $x\in M$, for which 
the rank of the sublattice spanned by the norm 4 vectors
in $\Gamma $ is strictly bigger than the one of the sublattice spanned
by the norm 4 vectors in $G(x)$, hence 
$\dim \langle \Gamma _4 \rangle \geq 7$.
\eb

{\bf Assumption:}
In view of Proposition \ref{a2a2} 
we will  assume in the following that all classes in 
${\cal K}$ are pairwise orthogonal. 

Then Lemma \ref{mink7} and Proposition \ref{classes} directly imply 

\begin{corollary}\label{llll}
$\Gamma $ has a sublattice $L\cong (A_1 \otimes A_2)^4$.
The lattice $\Lambda \cap \R L $ contains a sublattice $L'$ 
which is isometric to $(A_1 \otimes A_2^* )^4$ (under the same isometry).
\end{corollary}

\begin{corollary}\label{overlat}
There is an even overlattice  $\tilde{\Gamma } $  containing
$\Gamma ^{(e)}$ of index $9$.
\end{corollary}

\bew
Let $L' = \langle x_i,y_i | 1\leq i \leq 4\rangle$ with
$(x_i,y_i) = 2/3$, $(x_i,x_i) = (y_i,y_i) = 4/3$  for all $i$ and
all other scalar products are $0$
 be the lattice of Corollary \ref{llll}.
Then $L'$ contains an even sublattice $\langle v_1 := x_1+x_2+x_3, 
v_2:=x_2-x_3+x_4 \rangle$
such that $v_1$ and $v_2$ are linearly independent modulo $(L')^*$.
Therefore $\tilde{\Gamma } := \langle \Gamma ^{(e)} , v_1,v_2 \rangle $
is the desired lattice.
\eb

Note that $\det (\tilde{\Gamma } ) \leq 4 \cdot 2097 / 81 \leq 103.6 $
and if $\Gamma = \Gamma ^{(e)}$, then $\tilde{\Gamma }$ is 
of determinant $\leq 25$.
If $\Gamma \neq \Gamma ^{(e)}$ then we will take the 
subdirect product with $A_1^*$ to obtain an even lattice of
even determinant  $\leq 50$ in dimension 13. 
We therefore treat the two cases separately.

\begin{proposition}\label{CT}
Assume that $\Gamma = \Gamma ^{(e)}$.
Then $\Gamma \cong CT$.
\end{proposition}

\bew
We may assume that $\Lambda $ is generated by its minimal vectors.
Since these have norm $\frac{4}{3}$, the 
Sylow $p$-subgroup of the discriminant group
$\Lambda / \Gamma $ is generated by isotropic classes for all primes
$p\neq 3$.
In particular, if some prime $p \neq 3$ divides $\det (\Gamma )$,
then $\tilde{\Gamma } $ has an integral overlattice of 
index $p$. In particular $p^2$ divides $\det (\tilde{\Gamma })$.
By \cite[Table  15.4, p. 387]{SPLAG} 
there are no even 12-dimensional lattices of determinant 
1, 2, 3, or 6.
Therefore $\det (\tilde{\Gamma }) \neq p^2, 2p^2, 3p^2 $ or $6p^2$ 
for all primes $p\neq 3$.
Since $\det(\tilde{\Gamma }) \leq 25$ one gets 
$\det(\tilde{\Gamma }) = 9 $ or $16$.
In the first case $\tilde{\Gamma}$ is in the 
genus of $A_2\perp A_2 \perp E_8$.
The genus contains 3 isometry classes, for only one the
dual lattice has minimum $\geq 4/3$. This class is 
$E_6 \perp E_6$.
The automorphism group of $E_6\perp E_6$ has 20 orbits on the
sublattices of index 3, only one of which consists of lattices of 
which the dual has minimum $\geq 4/3$.
Continuing with this lattice, one finds 32 orbits 
of sublattices of index 3 under the automorphism 
group for only one of which the dual has minimum $4/3$.
This lattice is the Coxeter-Todd lattice $CT$.

In the second case, $\tilde{\Gamma }$ has an even overlattice 
$\tilde{\tilde{\Gamma }}$ which is in the genus of $E_8 \perp D_4$.
The genus of $E_8\perp D_4$ consists of 2 classes, the other is
$D_{12}$.
Both lattices have vectors of norm 1 in their duals, so this case 
is impossible.
\eb

Now assume that $[\Gamma : \Gamma ^{(e)} ] = 2$.
We then have the following situation:

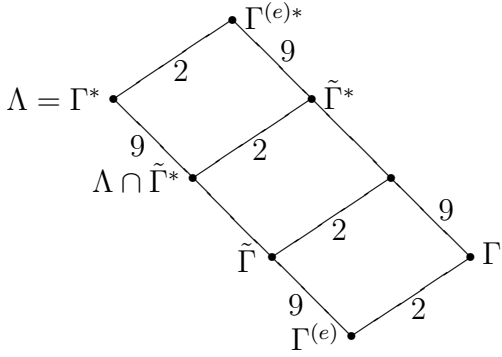
\setlength{\unitlength}{1.5pt}
\begin{picture}(120,110)(0,20)
\put(80,40){\circle*{2}} 
\put(95,47){\makebox(0,0)[l]{$2 $}}
\put(77,40){\makebox(0,0)[r]{$\Gamma ^{(e)} $}}
\put(110,60){\circle*{2}} 
\put(113,61){\makebox(0,0)[l]{$\Gamma $}}
\put(60,60){\circle*{2}} 
\put(68,48){\makebox(0,0)[r]{$9 $}}
\put(56,60){\makebox(0,0)[r]{$\tilde{\Gamma}  $}}
\put(75,67){\makebox(0,0)[l]{$2 $}}
\put(90,80){\circle*{2}} 
\put(40,80){\circle*{2}} 
\put(28,88){\makebox(0,0)[r]{$9 $}}
\put(37,80){\makebox(0,0)[r]{$\Lambda \cap \tilde{\Gamma} ^* $}}
\put(55,87){\makebox(0,0)[l]{$2 $}}
\put(70,100){\circle*{2}} 
\put(73,101){\makebox(0,0)[l]{$\tilde{\Gamma }^{*} $}}
\put(20,100){\circle*{2}} 
\put(17,100){\makebox(0,0)[r]{$\Lambda = \Gamma ^* $}}
\put(35,107){\makebox(0,0)[l]{$2 $}}
\put(50,120){\circle*{2}} 
\put(53,121){\makebox(0,0)[l]{$\Gamma ^{(e)*} $}}

\put(62,112){\makebox(0,0)[l]{$9 $}}
\put(102,72){\makebox(0,0)[l]{$9 $}}
\drawline(80,40)(20,100)(50,120)(110,60)(80,40)
\drawline(40,80)(70,100)
\drawline(60,60)(90,80)
\end{picture}

Let  $\varphi : \Gamma / \Gamma ^{(e)} \to A_1^*/A_1 $ be the unique
isomorphism. Then we define 
$$M:= \Gamma \subdir{} A_1^* := \{ x+y \in \Gamma \perp A_1^* \mid
\varphi (x) = y \}  .$$
The overlattice $\tilde{M} := \tilde{\Gamma } + M$ contains $M$ of
index 9 and is an even lattice in dimension 13 of even determinant
$\leq 50$. 
Moreover $M$ contains a unique pair $\pm v $ of vectors of norm 2, 
The orthogonal complement $v^{\perp }:= \{ \gamma \in M \mid (v,\gamma ) = 0 \}$
is the lattice $\Gamma ^{(e)}$.
The orthogonal complement of $v$ in $\tilde{M}$ is a sublattice of 
$\Lambda $ and hence has minimum $\geq \frac{4}{3}$. 

We will show that this situation is impossible:

\begin{proposition}
$\Gamma = \Gamma ^{(e)}$.
\end{proposition}

\bew
Assume that $\Gamma \neq \Gamma ^{(e)}$ and construct the lattice
$\tilde{M}$ as above.
Then $\det (\tilde{M}) = \frac{1}{2} \det(\tilde{\Gamma })$ is
an even integer $\leq 50$.

Moreover the condition (MIN) is satisfied:
\\
(MIN):
There is a vector $v\in \tilde{M}_2$ such that 
$L_v:=\{ \gamma \in M^* \mid (\gamma ,v ) = 0 \}$ is a lattice of
minimum $\geq \frac{4}{3}$. 

As in the proof of Proposition \ref{CT} 
we may assume that $\Lambda $ is generated by its minimal vectors.
Since these have norm $\frac{4}{3}$, the orthogonal complement $O$ of
$\Gamma +A_1 /\Gamma ^{(e)} + A_1 $ in $M^*/M$ is generated 
by classes of norm $\frac{4}{3}$.
In particular for all primes $p\neq 3$, the Sylow $p$-subgroup of  $O$
is generated by isotropic classes.
Hence if some prime $p \neq 3$ divides $\frac{1}{2} \det (\tilde{M})$,
then $\tilde{M} $ has an integral overlattice of 
index $p$. In particular $p^2$ divides $\frac{1}{2} \det (\tilde{M})$.
By \cite[Table  15.4, p. 387]{SPLAG} 
there are no even 13-dimensional lattices of determinant  2.
Therefore $\det (\tilde{M}) \neq 2p^2$ 
for all primes $p\neq 3$.
Since $\det(\tilde{M }) \leq 50$ one gets 
$\det(\tilde{M}) = 6 $, $16$, $18$, $24$, $32$ or $48$.

{\bf det=6}: 
In the first case $\tilde{M}$ is in the 
genus of $A_5 \perp E_8$.
The genus contains 3 isometry classes two of which, say $L_1$ and $L_2$
(namely all but $A_5\perp E_8$) satisfy condition (MIN).
No sublattice of index 9 of $L_1$ satisfies (MIN) and a unique
sublattice of index 9, say $L'$, satisfies the condition (MIN).
$L'$ has root system $A_1^3$, hence $L' \neq M$.

{\bf det=24}:
If $\det(\tilde{M}) = 24$, then $M$ is a sublattice of 
index 2 of the lattice $L'$ constructed in the case that 
$\det(M)  = 6$. But no such sublattice satisfies condition (MIN).

{\bf det=16}:
Assume now that $\det(\tilde{M})  = 16$.
Then $\tilde{M} $ has an even overlattice $\tilde{\tilde{M}}$ containing
$M$ of index $2$. 
The determinant of $\tilde{\tilde{M}}$ is 4.
By \cite[Table 15.4]{SPLAG} 
there is a unique genus of even 13-dimensional lattices of deterimant 4
 the genus of $D_{13}$. It contains 3 classes, none of which satisfies 
the condition (MIN).

{\bf det=18}:
Now assume that $\det{\tilde{M}} = 18$.
Then $\tilde{M} $ is a maximal integral lattice (since there is no
even lattice of determinant 2 in dimension 13).
In particular the discriminant group of $\tilde{M}$ is not cyclic, and
the 3-Sylow subgroup of $\tilde{M}^*/\tilde{M}$ is isometric to 
the unique anisotropic quadratic space of dimension 2 over $\F_3 $.
There is a unique genus of such lattices, namely the one of
$E_6 \perp E_6 \perp A_1$. 
This genus has class number 7. Only 
the lattice $E_6 \perp E_6 \perp A_1$  satisfies condition (MIN). 
Only one sublattice of index 9 of this lattice satisfies 
condition (MIN) and this lattice is isometric to $A_1 \perp CT $.

{\bf det=32}:
If $\det{\tilde{M}} = 32 $, then $\tilde{M}$ has an even overlattice
$\tilde{\tilde{M}} $ of determinant $8$. 
By \cite[Table 15.4]{SPLAG} there is a unique genus of even lattices
of determinant 8 in dimension 13, the one of $A_1\perp D_4 \perp E_8$.
This genus contains 4 classes, none of which satisfies condition (MIN).

{\bf det=48}:
Assume finally that $\det(\tilde{M}) = 48$.
Then there is an even overlattice $\tilde{\tilde{M}}$ of $\tilde{M}$ 
of determinant $\det(\tilde{\tilde{M}}) = 12$.
There is one genus of such lattices of determinant 12, namely the
one of $E_6\perp D_7$. 
Its class number is 7 and none of the lattices satisfies condition (MIN).
\eb

Together with Proposition \ref{CT} this concludes the proof of
Theorem \ref{s378}.

\renewcommand{\arraystretch}{1}
\renewcommand{\baselinestretch}{1}
\large

\begin{center}
Lehrstuhl D f\"ur Mathematik, RWTH Aachen,
Templergraben 64, 52062 Aachen, Germany,
\\
e-mail: nebe@math.rwth-aachen.de
\end{center}

\begin{center}
St. Petersburg Branch of the Steklov Mathematical Institute,
Fontanaka 27, 191011 St. Petersburg, Russia,
\\
e-mail: bbvenkov@yahoo.com
\end{center}

\end{document}